\documentclass[10pt]{amsart}
\usepackage{graphicx}
\usepackage{psfrag}
\usepackage{amsthm}
\usepackage{amscd}
\usepackage{amsfonts}
\usepackage{amstext}
\usepackage{amssymb}
\usepackage{plain}
\def\cal{\mathcal}

\let\newpf\proof \let\proof\relax 
\newenvironment{pf}{\newpf[\proofname]}{\qed\endtrivlist}

\def\0{{\mathbf{0}}}

\newcommand{\Dist}{\operatorname{Dist}}

\def\g{\gamma}
\def\tg{{\tilde\gamma}}
\def\a{{\tilde a}}
\def\b{{\tilde b}}

\def\d{{\underline d}}

\newcommand{\ntop}[2]{\genfrac{}{}{0pt}{1}{#1}{#2}}

\newtheorem*{Theorem A}{Theorem A}
\newtheorem*{Theorem B}{Theorem B}
\newtheorem*{Theorem C}{Theorem C}

\newtheorem{thm}{Theorem}[section]

\newtheorem{lemma}[thm]{Lemma}

\newtheorem{prop}[thm]{Proposition}

\theoremstyle{remark}
\newtheorem{rem}{Remark}[section]

\numberwithin{equation}{section}

\def \bn {\hfill \\ \smallskip\noindent}

\theoremstyle{definition}

\def\proof{\bn {\bf Proof.} }

\newcommand{\inter}{\operatorname{int}}

\newcommand{\id}{\operatorname{id}}

\newcommand{\AAA}{{\cal A}}

\newcommand{\UUU}{{\cal U}}
\newcommand{\VV}{{\cal V}}

\newcommand{\C}{{\mathbb C}}

\newcommand{\N}{{\mathbb N}}

\newcommand{\R}{{\mathbb R}}

\newcommand{\Z}{{\mathbb Z}}

\def\B0{{\bold{0}}}

\catcode`\@=12

\def\Empty{}
\newcommand\oplabel[1]{
  \def\OpArg{#1} \ifx \OpArg\Empty {} \else
        \label{#1}
  \fi}
                
\newcommand{\comm}[1]{}
\newcommand{\comment}[1]{}

\begin{document}

\title[Statistical properties in smooth families of unimodal maps]
{Statistical properties of unimodal maps: smooth families with
negative Schwarzian derivative}

\author{Artur Avila and Carlos Gustavo Moreira}

\address{
Coll\`ege de France, 3 Rue d'Ulm, 75005 Paris -- France.
}

\email{avila@impa.br}

\address{
IMPA -- Estr. D. Castorina 110 \\
22460-320 Rio de Janeiro -- Brasil.
}
\email{gugu@impa.br}

\thanks{Partially supported by Faperj and CNPq, Brazil.}

\date{\today}

\begin{abstract}

We prove that there is a residual set of families of
smooth or analytic unimodal maps with quadratic critical point
and negative Schwarzian derivative such that
almost every non-regular parameter is Collet-Eckmann with
subexponential recurrence of the critical orbit.  Those conditions lead to
a detailed and robust statistical description of the dynamics.
This proves the Palis conjecture in this setting.

\end{abstract}

\maketitle

\section{Introduction}

`The main strategy of the study of all mathematical
models is, according to Poincar\'e, the consideration of each model as a
point of a space of different but similar admissible systems'(V. Arnold in
\cite {Ar}).  One of the main concerns of dynamical systems is to
establish properties valid for typical systems.  Since the space of
such systems is usually infinite dimension, there are of course many
concepts of `typical'. 
According to \cite {Ar} again, `The most physical genericity notion   
is defined by Kolmogorov (1954), who suggested to call a property of
dynamical systems exceptional, if it holds only on Lebesgue measure
zero set of values of the parameters in every (topologically) generic
family of systems, depending on sufficiently many parameters'.

In the last decade Palis \cite {Pa} described a general program for
(dissipative) dynamical systems in any dimension.
He conjectured that a typical dynamical
system has a finite number of attractors described by physical
measures, the union of their basins has full Lebesgue measure, and
those physical measures are stochastically stable.  Typical was to be
interpreted in the Kolmogorov sense: full measure in generic families.
Our aim here is to give a proof of this conjecture
for an important class of one-dimensional dynamical systems.

Here we consider unimodal maps, that is, continuous maps from an interval to
itself which have a unique turning point.  More specifically, we consider
$S$-unimodal maps, that is, we assume that the map is $C^3$
with negative Schwarzian derivative and
that the critical point is non-degenerate.

\subsection{The quadratic family}

The basic model for unimodal maps is the quadratic family, $q_a(x)=a-x^2$,
where $-1/4 \leq a \leq 2$ is a parameter.  Despite its simple appearance,
the dynamics of those maps presents many remarkable phenomena.  Restricting
to the probabilistic point of view, its richness first became apparent with
the work of Jakobson \cite {J},
where it was shown that a positive measure set of parameters
corresponds to quadratic maps with stochastic behavior.  More precisely,
those parameters possess an absolutely continuous invariant
measure (the physical measure) with positive Lyapunov exponent.
On the other hand, it was later shown by
Lyubich \cite {puzzle} and Graczyk-Swiatek \cite {GS} that regular
parameters (with a periodic hyperbolic attractor) are (open and) dense.
So at least two kinds of very distinct observable behavior are present
on the quadratic family, and they alternate in a complicate way.

Besides regular and stochastic behavior, different behavior was shown to
exist, including examples with bad statistics, like absence of a
physical measure or a physical measure concentrated on a
hyperbolic repeller.  Those pathologies were shown to be
non-observable in \cite {parapuzzle} and \cite {MN}.  Finally in \cite
{regular} it was proved that almost every real quadratic map is either
regular or stochastic.

Among stochastic maps, a specific class grabbed lots of attention in the
90's: Collet-Eckmann maps.  They are characterized by a positive Lyapunov
exponent for the critical value, and gradually they were shown to have `best
possible' near hyperbolic properties: exponential decay of correlations,
validity of central limit and large deviations theorems, good spectral
properties and zeta functions (\cite {KN}, \cite {Y2}).  Let us call
attention to the robustness of the statistical description, with a good
understanding of stochastic perturbations: strong stochastic stability
(\cite {BV}), rates of convergence to equilibrium (\cite {BBM}).

In \cite {AM} the regular or stochastic dichotomy was extended by showing
that almost every
stochastic map is actually Collet-Eckmann and has polynomial
recurrence of its critical point,
in particular implying the validity of the above mentioned results.

The position of the quadratic family in the borderline of real and complex
dynamics made it a meeting point of many different techniques:
most of the deeper results depend on this interaction.  It gradually became
clear however that studying the quadratic family allows one
to obtain results on
more general unimodal maps.

\subsection{Universality}

Starting with the works of Milnor-Thurston, and also through the discoveries
of Feigenbaum and Coullet-Tresser, the quadratic family was shown to be
a prototype for other families of unimodal maps which presents universal
combinatorial and geometric features.  More recently, the result of density
of hyperbolicity among unimodal maps was obtained in \cite {K} exploiting
the validity of this result for quadratic maps.

In \cite {ALM}, a general method was developed
to transfer information from
the quadratic family to real analytic families of unimodal maps.  It was
shown that the decomposition of spaces of analytic unimodal maps according
to combinatorial behavior is essentially a codimension-one lamination.    

Thinking of two analytic families as transversals to this lamination, one
may try to compare the parameter space of both families via the holonomy 
map.  A straightforward application of this method allows one to conclude   
that the bifurcation pattern of a general analytic family is locally     
the same as in the quadratic family from the
topological point of view (outside of countably many `bad parameters').

The `holonomy' method was then successfully applied to
extend the regular or stochastic dichotomy from the quadratic family to a
general analytic family.  The probabilistic point of view presents new   
difficulties however.  First, the statistical properties of two
topologically conjugate maps need not correspond by the
(generally not absolutely continuous) conjugacy.  Fortunately many
properties are preserved, in
particular the criteria used by Lyubich in his result.

The second difficulty is that the holonomy map is usually not absolutely
continuous, so typical combinatorics for the quadratic family may not
be typical for other families: it has to be shown
that the class of   
regular or stochastic maps is still typical after application of the
holonomy map.

\subsection{Results and outline of the proof}

Let us call a $k$-parameter family good if almost every non-regular
parameter is Collet-Eckmann (and satisfies
some additional technical conditions).
Our goal will be to prove that good families are generic.  This question
naturally makes sense in different spaces of unimodal maps (corresponding to
different degrees of smoothness).  We only deal with the last steps of this
problem (going from the quadratic family to analytic and then smooth
categories), basing ourselves on the building blocks \cite
{parapuzzle}, \cite {regular}, \cite {ALM}, and \cite {AM}.

We start by describing how the holonomy method of \cite {ALM}
can be applied to generalize the results of \cite {AM} to
general analytic families (to put together those two papers we need to
do a non-trivial strengthening of \cite {AM}).
As a consequence we conclude that essentially
all analytic families are good.

To get to the smooth setting (at least $C^3$, since we are assuming negative
Schwarzian derivative), our strategy is different: we show a certain
robustness of good families, which together with their denseness (due to the
analytic case) will yield genericity.  Our main tool is one of the nice
properties of Collet-Eckmann maps: persistence of the Collet-Eckmann
condition under generic unfolding (a result of \cite {T1}).  By means of
some general argument, we reduce the global result to this local one.

Let us mention that the results of this paper
are still valid without the negative
Schwarzian derivative assumption (also allowing one to get to
$C^2$ smoothness), see \cite {A}, \cite {AM4}.
The techniques are very different however, since we replace the
global holonomy method we use here by a local
holonomy analysis based on a ``macroscopic'' version of the
infinitesimal perturbation method of \cite {ALM}.  For analytic
maps this also allowed us to obtain better asymptotic
estimates which have interesting consequences, for instance pathological
measure-theoretical behavior of the lamination by combinatorial classes (see
\cite {AM2}).

{\bf Acknowledgements:}
We thank Viviane Baladi, Mikhail Lyubich, and Marcelo Viana for helpful
discussions and suggestions.

\section{General definitions}

\subsection{Notation}

Let $I=[-1,1]$ and let
$B^k$ be the closed unit ball in $\R^k$ (we will use the notation $I$ for
the {\it dynamical} interval, while $B^1$ will be reserved for the
one-dimensional {\it parameter} space).  We will consider $B^k$
endowed with the Lebesgue measure normalized so that $|B^k|=1$.
Let $C^r(I)$ denote the space of $C^r$ maps $f:I \to \R$.

By a {\bf unimodal map} we will mean a smooth (at least $C^2$)
symmetric (even)
map $f:I \to I$ with a unique critical point at $0$ such that $f(-1)=-1$,
$Df(-1) \geq 1$, and if $Df(-1)=1$ then $D^2 f(-1)<0$.
If $f$ is $C^3$, we define the Schwarzian derivative on $I
\setminus \{0\}$ as
$$
Sf=\frac {D^3 f} {Df}-\frac {3} {2} \left ( \frac {D^2 f} {Df} \right )^2.
$$

For $a>0$, let $\Omega_a \subset \C$ denote an $a$ neighborhood $I$.

Let $\AAA_a$ denote the space of holomorphic maps on
$\Omega_a$ which have a continuous extension to $\partial \Omega_a$,
satisfying $\phi(z)=\phi(-z)$, $\phi(-1)=\phi(1)=-1$ and $\phi'(0)=0$.

Notice that $\AAA_a$ is a closed affine subspace of the Banach space
of bounded holomorphic maps of $\Omega_a$.  We endow it with the induced
metric and affine structure.

We define $\AAA^\R_a \subset \AAA_a$ the space of maps which are real
symmetric.

\subsection{More on unimodal maps}

A $C^3$ unimodal map such that $Sf<0$ on $I \setminus \{0\}$ and such that
its critical point is non-degenerate (that is, $D^2 f \neq 0$)
will be called a {\bf $S$-unimodal map}.

We say that $x$ is a periodic orbit (of period $n$) for $f$ if $f^n(x)=x$
and $n \geq 1$ is minimal with this property.  In this case we define
$Df^n(x)$ as the multiplier of $x$.  Notice that this definition
depends only on the orbit of $x$.  We say that $x$ is
hyperbolic if $|Df^n(x)| \neq 1$.

A unimodal map is called {\bf regular} (or hyperbolic) if it has a
hyperbolic periodic attractor.  This condition is $C^2$ open.

A {\bf $k$-parameter family} of unimodal maps is a map
$F:B^k \times I \to I$ such
that for $p \in B^k$, $f_p(x)=F(p,x)$ is a unimodal map.
Such a family is said to be $C^n$ or analytic, according to $F$ being  
$C^n$ or analytic.  We introduce the natural topology in spaces of smooth   
families ($C^n$ with $n=2,...,\infty$), but do not introduce any topology   
in the space of analytic families (however, we will refer from time to time
to induced $C^n$ topologies).

An analytic family of $S$-unimodal maps $F$ will be called {\bf non-trivial}
if there exists a regular parameter.  Notice that this condition is
$C^3$-open.

A unimodal map $f$ is called {\bf Collet-Eckmann}
(CE) if there exists constants
$C>0$, $\lambda>1$ such that for every $n>0$,
$$
|Df^n(f(0))|>C \lambda^n.
$$
This means that the map is strongly hyperbolic along the critical orbit.
It is also useful to study the hyperbolicity of 
backward iterates of the critical point, so we say that $f$ is {\bf
Backwards Collet-Eckmann} (BCE) if there exists $C>0$, $\lambda>1$ such
that for any $n>0$ and any $x$ with $f^n(x)=0$, we have
$$
|Df^n(x)|>C \lambda^n.
$$
By a result of Nowicki (see \cite {MvS}), for $S$-unimodal maps CE implies
BCE, so we will mostly discuss the Collet-Eckmann condition (except for the
last section where we consider $C^2$ unimodal maps as well).

Very often it is useful to estimate how fast is the recurrence of the
critical orbit.  We will be mainly interested in two kinds of control:
{\bf Polynomial Recurrence} (P) if there exists $\alpha>0$ such that
$$
|f^n(0)|>n^{-\alpha}
$$
for big enough $n$ and {\bf Subexponential Recurrence} (SE) if for all
$\alpha>0$,
$$
|f^n(0)|>e^{-\alpha n}
$$
for $n$ big enough.

We will say that $f$ is {\bf Weakly Regular} (WR) if
$$
\lim_{\delta \to 0} \liminf_{n \to \infty} \frac {1} {n}     
\sum_{\ntop {1 \leq k \leq n} {f^k(0) \in (-\delta,\delta)}}
\ln |Df(f^k(0))|=0.
$$
This condition is used in proofs of
stochastic stability for $C^2$ maps, see \cite {T2}.

We will consider spaces of $S$-unimodal maps: we define
$\UUU^r \subset C^r(I)$ the set of
$S$-unimodal maps.  Spaces of analytic unimodal maps are now easily defined:
$\UUU_a=\UUU^3 \cap \AAA^\R_a$.

\subsection{The quadratic family}

The quadratic family is the most studied family of unimodal maps.  It is
usually parametrized by
$$
q_t(x)=t-x^2,
$$
so that for $-1/4 \leq t \leq 2$, there exists a unique symmetric
interval $I_t=[-\beta_t,\beta_t]$ such that
$q_t(I_t) \subset I_t$ and $q_t(-\beta_t)=-\beta_t$,
so $q_t$ can be seen as a unimodal map of $I_t$ (which depends on $t$).
Moreover $Sq_t(x)<0$ if $x \neq 0$.

By an affine reparametrization of the parameter $t$ and of each interval
$I_t$, we obtain a canonical one-parameter family of $S$-unimodal maps in
the interval $I$, which we denote $p_t$, $t \in B^1$, which will be called
the quadratic family as well.

\subsection{Quasisymmetric maps}

Let $\g \geq 1$ be given.
We say that a homeomorphism $f:\R \to \R$ is {\bf quasisymmetric} (qs)
if there exists a constant $k>1$ such that for all $x \in \R$ and
any
$h>0$
$$
\frac {1} {k} \leq \frac {f(x+h)-f(x)} {f(x)-f(x-h)} \leq k.
$$

A homeomorphism $h$ is quasisymmetric if and only if it admits a
real-symmetric extension to a quasiconformal map $\tilde h:\C \to \C$
(Ahlfors-Beurling).
We will say that $h$ is $\g$-qs (or that $\g$ is a qs constant for $h$) if
the dilatation of $\tilde h$ is bounded by $\g$.  This definition of the
quasisymmetric constant is convenient since the composition of
quasisymmetric maps $g$ and $f$ is readily seen
to be quasisymmetric and the qs constant of $g \circ f$ is bounded by the
product of the qs constants of $g$ and $f$.

If $X \subset \R$ and $h:X \to \R$ has a $\g$-quasisymmetric extension to
$\R$ we will also say that $h$ is $\g$-qs.

\section{Statement of the results}

\subsection{A dichotomy for generic families of $S$-unimodal maps}

We would like to classify the typical behavior in generic
families of unimodal maps.  This classification should reveal refined
information on the stochastic description of the dynamics of those typical
parameters.

We will therefore consider a smooth enough family of unimodal maps $F$.
The techniques of the present paper will need the fact that
$F$ is a family of $S$-unimodal maps.  This includes two main
restrictions: the negative Schwarzian derivative
and the quadratic critical point.
The first one is serious, since this condition is not dense, but can be
removed with more refined techniques (see \cite {A}).  The second one
(which is not present in the usual definition of $S$-unimodal map, but is
rather a convention in this paper) is no serious loss of
generality, since quadratic critical point is certainly typical
among unimodal maps.

\begin{rem}

Families of unimodal maps with a fixed critical exponent different from $2$
have also been subject of much study.  This theory has many similarities,
but also some important differences and new features,
and is not nearly as complete as the case of criticality $2$.  It is however
widely expected that the Palis conjecture
(and indeed our Theorems A, B and C) still holds in this setting.

\end{rem}

We first consider the analytic case.

\begin{Theorem A}

Let $F$ be a non-trivial $k$-parameter
analytic family of $S$-unimodal maps.
Then for almost every non-regular parameter
$p \in B^k$, $f_p$ satisfies the Collet-Eckmann and
Polynomial Recurrence conditions.

\end{Theorem A}

Notice that the set of non-trivial analytic
families is indeed generic in any meaningful sense: its
complement has ``infinite codimension'', see Proposition \ref
{strongdense}.  Moreover, if an analytic
family is non-trivial, it is possible to verify the non-triviality
in finite time (with an infinite precision computer\footnote {Since regular
parameters form an open set (non-empty if the family is non-trivial), and
any regular parameter one can be also checked in finite time
(by locating the attracting hyperbolic periodic orbit.}).

\comm{
In this case, the generic families we
will consider will be
families satisfying a non-degeneracy condition, that is, we assume
they contain a regular map (non-trivial families).
This condition is easy to verify, both in practice (if the
family is non-degenerate this can be verified in finite time) as in theory
(they are indeed generic for meaningful concepts, see Proposition \ref
{strongdense}).  Our main result is a dichotomy on the
typical dynamical behavior in any such family.
}

Our second result about non-trivial analytic
families is the robustness of a slightly
weaker dichotomy under $C^2$ perturbations of the family.

\begin{Theorem B}

Let $F$ be a non-trivial $k$-parameter analytic
family of $S$-unimodal maps.  Let
$F^{(n)}$ be a sequence of $C^2$ families such that
$F^{(n)} \to F$ in the $C^2$ topology.  For each $n$, let
$X_n$ be the set of parameters
$p \in B^k$ where $F^{(n)}$ is either regular or has only repelling periodic
orbits and satisfies simultaneously
the Backwards Collet-Eckmann, Collet-Eckmann, Subexponential
Recurrence and Weak Regularity conditions.  Then $|X_n| \to 1$.
In particular, almost every parameter of $F$ is Weakly Regular.

\end{Theorem B}

As a consequence, we can use a Baire argument to conclude that the dichotomy
is still valid among topologically generic smooth families (that is,
belonging to some residual set),
obtaining the following corollary of Theorems A and B.

\begin{Theorem C} [Smooth Dichotomy]

In topologically generic $k$-parameter $C^r$, $r=3,4,...,\infty$
families of $S$-unimodal maps,
almost every non-regular parameter satisfies the Backwards
Collet-Eckmann, Collet-Eckmann, Subexponential Recurrence and Weak
Regularity conditions.

\end{Theorem C}

It is good to recall that both types of behavior described by the
dichotomy are indeed observable for open sets of families of unimodal maps
(\cite {J}, \cite {BC}).

\begin{rem}

The space of $S$-unimodal maps is easy to describe and easier to work with
but has some disadvantages.  One of them is that it is not an intrinsic
condition, in particular it is not invariant by analytic change of
coordinates.  A more natural class to work with is the space of
quasiquadratic unimodal maps as defined by \cite {ALM}.  A unimodal map
$f$ is called quasiquadratic if there exists a $C^3$-neighborhood of $f$
where all maps are topologically conjugate to some quadratic map.
The results of this paper are still valid in spaces of quasiquadratic
unimodal maps (which includes $S$-unimodal maps).  The proofs are unchanged,
since the results we need from \cite {ALM} are stated
and proved for quasiquadratic maps.  We remark further that the
description of quasiquadratic unimodal maps can be used to describe all
unimodal maps: it is proved in \cite {A}, \cite {AM4}
that (Kolmogorov) typical (analytic or smooth) unimodal maps
have either a quasiquadratic renormalization or a
quasiquadratic unimodal restriction.

\end{rem}

\subsection{Ergodic consequences}

The importance of the above dichotomy is the fact that each of the two
possibilities has very well defined stochastic properties.  We quickly
recall those (we assume that maps are $S$-unimodal).

Regular maps have a periodic attractor whose basin is big both topologically
(open and dense set) as in the measure-theoretical sense (full measure). 
Moreover the attractor and its basin are stable under $C^1$ perturbations. 
The dynamics of such maps can be described in deterministic terms.

Maps satisfying CE and SE have non-deterministic dynamics.  They can be
however described
through their stochastic properties, and it turns out that
such maps have the main good properties usually found in hyperbolic maps.
First, there is a physical measure, that is an invariant probability which
describes asymptotic behavior of orbits: for almost every $x$ and for every
continuous $\phi:I \to \R$,
$$
\lim \frac {1} {n} \sum_{k=0}^{n-1} \phi(f^k(x))=\int \phi d\mu.
$$
This physical measure has a positive Lyapunov exponent and
is indeed absolutely continuous and supported on a
cycle of intervals, so the asymptotic behavior is non-deterministic.
The convergence to the asymptotic stochastic model is exponential, see the
results on decay of correlations and convergence to equilibrium
(\cite {KN}, \cite {Y2}).  Those properties are beautifully related to a
spectral gap of a transfer operator and to zeta functions, see \cite {KN}. 
Notice finally that exponential decay of correlations is actually equivalent
to the Collet-Eckmann condition (see \cite {NS}).

While the dynamics is highly unstable under deterministic perturbations
(nearby maps can be regular for instance), the stochastic
description given by the physical measure $\mu$ is robust under
stochastic perturbations: the
perturbed system has a stationary measure which is close to $\mu$ in the
sense of the $L^1$ distance between their densities (\cite {BV}).
For studies of decay
of correlations for the perturbed systems, see \cite {BBM}.

\section{Analytic families}

\subsection{Hybrid classes and holonomy maps}

Two $S$-unimodal maps $f, \tilde f$
are said to be {\bf hybrid equivalent}
if they are topologically conjugate and, in case they are regular, their
attracting periodic orbits have the same multiplier.

The set of all maps which are hybrid equivalent to some $f$ is called
the {\bf hybrid class} of $f$.  The partition of $S$-unimodal maps
into hybrid classes is thus a refinement of the partition in topological
conjugacy classes.

It follows from a result of Guckenheimer (see \cite {MvS})
that any $S$-unimodal map $f$ is topologically
conjugate to some quadratic map.  It turns out that if $f$ has a hyperbolic
attractor, we can select the quadratic map with a hyperbolic attractor with
the same multiplier\footnote {This follows for instance from
Milnor-Thurston kneading theory and the fact that the quadratic family is a
full family.  Another way to see this is to notice that
in each ``hyperbolic window'' of quadratic maps
(a maximal parameter interval $(a,b)$ such that $p_t$ is hyperbolic
for $t \in (a,b)$), the multiplier of the hyperbolic attractor induces a
homeomorphism from $(a,b)$ to $(-1,1)$ (this is a consequence for instance
of the work of Douady-Hubbard on the complex quadratic family).}.
In particular, each hybrid class intersects the
quadratic family in at least one point.

The problem of uniqueness is much harder.
The following result is due to Lyubich \cite {puzzle}
and Graczyk-Swiatek \cite {GS}, and is a consequence of (the proof of)
the equivalent rigidity result for quadratic maps:

\begin{thm} \label {phase}

Let $h$ be a topological
conjugacy between two analytic $S$-unimodal maps
$f$ and $\tilde f$ which have all periodic orbits repelling.
Then $h$ is quasisymmetric.

\end{thm}

\begin{rem}

Although we won't use it here, a similar theorem still holds for maps with
non-repelling periodic orbits: if $f$ and $\tilde f$ are two topologically
conjugate $S$-unimodal maps
and have non-repelling periodic orbits then we can select a
topological conjugacy which is quasisymmetric (the choice of the topological
conjugacy is not unique).  This result is considerably easier
than the case where all periodic orbits are repelling, and does not use
analyticity.

\end{rem}

This rigidity result has a remarkable consequence for quadratic maps: each
hybrid class intersects the quadratic family at a unique parameter.
Thus, any $S$-unimodal map $f$ is hybrid equivalent
to a unique quadratic map $\chi(f)$.  The map $\chi$ is called
the {\bf straightening}
\footnote {We should point out that there is also a notion of hybrid class
in complex dynamics.  In that context, the fact that each hybrid class (of
quadratic-like maps with connected Julia set) contains exactly one quadratic
polynomial is a consequence of the Straightening Theorem of Douady-Hubbard.
Our definition of hybrid class is motivated precisely by the possibility of
defining an analogous straightening map (whose existence is proved by quite
different methods).}.

\begin{lemma} \label {qs}

Let $f$ be an analytic $S$-unimodal map.  Then $\chi(f)$ is regular/CE/P if
and only if $f$ also satisfies the corresponding property.

\end{lemma}

\begin{pf}

The property of being regular is clearly invariant under hybrid equivalence,
so we only have to analyze invariance of the conditions CE and P.

By \cite {NP}, the Collet-Eckmann condition is topologically
invariant, so it is preserved under hybrid equivalence.

To check invariance of polynomial recurrence of the critical orbit,
first assume that $f$ has some non-repelling periodic point $p$.
In this case, the the orbit of $p$ must attract the
critical point.  In particular, the critical point is either non-recurrent
(in which case both $f$ and $\chi(f)$ satisfy P in a trivial way) or
periodic (in which case $f$ and $\chi(f)$ do not satisfy P also in a trivial
way).

If $f$ has all periodic
orbits repelling, by Theorem \ref {phase}, the conjugacy between $f$ and
$\chi(f)$ is quasisymmetric, and in particular H\"older.  It is easy to see
that P is invariant by H\"older conjugacy.
\end{pf}

\begin{rem}
By \cite {NP1}, two $S$-unimodal Collet-Eckmann maps which are
topologically conjugate are H\"older conjugate, so using \cite {NP} we see
that the joint conditions CE and P
are topologically invariant.
This joint invariance of CE and P is all that will
be used in the further arguments.  Notice that \cite {NP1} and \cite
{NP} do not assume analyticity, and are more elementary
than Theorem \ref {phase}.
\end{rem}

\subsection{Hybrid laminations}

It is natural to study the hybrid class of some map $f$.  This is what is
done in Theorem A of \cite {ALM} in the analytic setting, where it is shown
that in $\UUU_a$, every hybrid class is a codimension-one analytic
submanifold.  Moreover, different hybrid class fit together in some nice
structure, called {\bf hybrid lamination}.

\begin{rem}

It is not known if the
hybrid lamination is really a lamination everywhere.  In \cite {ALM},
it is shown that the
hybrid lamination is a lamination (in the usual sense) ``almost
everywhere'' (more precisely, if restricted to an open set containing the
complement of countably many classes corresponding to existence of neutral
periodic orbits), which is enough for our purposes.

\end{rem}

A $k$-parameter analytic family of $S$-unimodal maps can be thought as an
analytic map from $B^k$ to some $\UUU_a$.
As a consequence, the structure of the hybrid lamination implies
that non-trivial analytic families are indeed quite frequent.

\begin{lemma} [Most analytic families are non-trivial] \label {strongdense}

If a $k$-parameter
analytic family of $S$-unimodal maps is not contained in some
non-regular hybrid class then it is non-trivial.
In particular, non-trivial analytic families are dense
in the space of $C^n$ families of $S$-unimodal maps, $n=3,...,\infty$.

\end{lemma}

\begin{pf}

Let us consider an analytic family of $S$-unimodal maps $F$.
By the theory of Milnor-Thurston, see \cite {MvS},
either all parameters have the
same non-periodic {\it kneading sequence}, or there exists a parameter with
periodic critical point.  In the latter case, the family is of course
non-trivial, so let us consider the former case.  Two $S$-unimodal maps with
the same kneading sequence are either topologically conjugate, or one of
them possess a neutral periodic orbit (see Corollary, Chapter 2,
page 157 of \cite {MvS}), and it follows that
the other is necessarily regular.  Thus,
if the family $F$ does not have regular parameters, all maps are non-regular
and topologically conjugate, that is, $F$ is contained in a non-regular
hybrid class.


For the denseness result, given a $C^r$ family $F$,
approximate it by an analytic family $\tilde F$.
If such an analytic family is contained in a hybrid class,
we can perturb it further in order to intersect two hybrid classes, since
each hybrid class is a codimension-one submanifold.
\end{pf}

Let us consider the case where $F$ is a one-parameter
analytic family of $S$-unimodal maps, that is,
an analytic curve in some $\UUU_a$.
A consequence of the nice structure of the hybrid lamination is the
following result:

\begin{lemma} [see the proof of Theorem C of \cite {ALM}] \label {analytic}

If $F$ is a one-parameter analytic family of $S$-unimodal maps which is not
contained in some hybrid class then there is an open set 
of parameters, with countable complement,
where $F$ is transverse to the hybrid lamination.

\end{lemma}

Define the map $\chi_F$ on $B^1$ by $\chi_F(t)=\chi(f_t)$.
In \cite {ALM} the map $\chi_F$ is considered as the holonomy map from
$F$ to the quadratic family along the hybrid lamination in some
$\UUU_a$.  Using this interpretation, they obtain the following result:

\begin{thm} [Theorem C of \cite {ALM}] \label {holonomy}

Let $F$ be a one-parameter family of unimodal maps which is not
contained in some hybrid class.  Then there is an
open set $U \subset B^1$ with countable complement such that the
straightening $\chi_F$ is quasisymmetric in any compact interval
$J \subset U$.

\end{thm}

\subsection{Dichotomy in the quadratic family}

The main result of \cite {AM} is that almost every parameter in the
quadratic family is either regular or Collet-Eckmann with a polynomial
recurrence of the critical orbit.  To obtain the
same result for a non-trivial analytic
family using Theorem \ref {holonomy},
we will need a stronger estimate, since quasisymmetric maps are not
in general absolutely continuous.

Let us say that a set $X \subset B^1$ has {\bf total qs-probability}
if the image
of $B^1 \setminus X$ by any quasisymmetric map $h:B^1 \to B^1$ has
zero Lebesgue measure.

By an improvement of the proofs in \cite {AM} (see appendix),
it is possible to obtain the following result:

\begin{thm} \label {capacity}

The set of quadratic maps which are either regular or
simultaneously CE and P has total qs-probability.

\end{thm}

\begin{rem}

In \cite {AM} a better result than polynomial recurrence
is obtained in the quadratic family.
Namely it is shown that the asymptotic exponent
of the recurrence
$$
\limsup_{n \to \infty} \frac {-\ln |f^n(0)|} {\ln n}
$$                           
is exactly $1$ for almost every non-regular map.  However, for a set of
total qs-probability, we are only able to show that the asymptotic exponent
is bounded.

\end{rem}

\subsection{Proof of Theorem A}

Let $F$ be a non-trivial analytic
family.  If all parameters are regular, there is
nothing to prove, so assume that there is a non-regular parameter.

First assume $F$ is one-parameter. 
By Theorems \ref {capacity} and \ref {holonomy}, for almost every $t \in
B^1$, $\chi_F(t)$ is either regular or satisfies CE and P.  By Lemma
\ref {qs}, this implies that $f_t$ is either regular or CE and P.

Assume now that $F$ is a $k$-parameter family.  Let $p \in B^k$ be a
regular parameter.  Let $L:B^1 \to B^k$ be an affine map such $p \in
L(B^1)$.  Let $F^L$ be the one-parameter family defined by
$f^L_t=f_{L(t)}$.  Then $F^L$ is a non-trivial one-parameter
analytic family and hence for almost every $t$, $f^L_t$ is
either regular or CE and P.  The result follows by application of Fubini's
Theorem.

\section{Robustness of the dichotomy}

To obtain the robustness claimed on Theorem B our approach
will be to exploit an important
result of Tsujii, whose core is a strong generalization of
Benedicks-Carleson result and techniques.  This result establishes that the
CE and SE conditions are infinitesimally persistent
in one-parameter families unfolding generically: they are density points of
CE and SE parameters.  The connection with our robustness
result, which has a global nature, is done using
some general argument.

\subsection{Tsujii's theorem}

Let $F$ be a $C^2$ $k$-parameter family of unimodal maps.  Assume that
$p_0$ is a parameter such that $f_{p_0}$ satisfies CE, BCE, SE, has
a quadratic critical point and all periodic orbits repelling. 
Tsujii's Theorem considers the case where $F$ is a generic unfolding
at $p_0$.  For one-parameter families, generic unfolding means precisely
\begin{equation} \label {eq}
\sum_{j=0}^\infty \frac {v(f_{p_0}^j(0))}
{Df_{p_0}^j(f_{p_0}(0))} \neq 0, \quad \text {where} \quad
v= \left. \frac {d} {dp} f_p \right |_{p=p_0}.
\end{equation}
This transversality condition will be called {\bf Tsujii transversality}.

If $F$ is a one-parameter family, we will say that $(F,p_0)$
satisfies the {\bf Tsujii conditions}
if all above requirements are satisfied.

The following is an immediate consequence of the main theorem of
Tsujii in \cite {T1}.

\begin{thm}

Let $F$ be a $C^2$ one-parameter family of unimodal maps.  Assume
$(F,t_0)$ satisfies the Tsujii conditions.  Then $t_0$ is a
density point of parameters $t$ for which $(F,t)$ satisfies the Tsujii
conditions and for which $f_t$ is WR.

\end{thm}

\subsection{A higher dimensional version}

In order to pass from one-parameter to $k$-parameters, we will
need the following easy proposition.
Let us say that $p \in B^k$ is a density point
of a set $X$ along a line $l$ through
$p$ if $p$ is a density point of $l \cap X$
in $l$ (endowed with the linear Lebesgue measure).

\begin{prop} \label {dimension}

If $p \in B^k$
is a density point of $X$ along almost every line, then $p$ is a
density point of $X$ in $B^k$.

\end{prop}

\begin{pf}

Let $E$ be the characteristic function of $X$.
For each line $l$ through $p$, let $A_l:\R \to l$ be an isometric
parametrization of $l$ taking $0$ into $p$.  Let $P^{k-1}$
be the space of such lines with the natural probability measure
(obtained by
identification with the $k-1$ dimensional projective space).  Let
$$
\rho_{\epsilon}(l)=\int_{-1}^{1} |r| E(A_l(\epsilon r)) dr.
$$
Assuming that $p$ is a density point of $X$ along almost every $l$ we have,
for almost every $l$
$$
\lim_{\epsilon \to 0} \rho_{\epsilon}(l)=1.
$$
Using polar coordinates, the relative measure of
$X$ in an $\epsilon$ ball around $p$ is given by
$$
\int_{P^{k-1}} \rho_{\epsilon}(l) dl.
$$
By the Lebesgue Convergence Theorem,
$$
\lim_{\epsilon \to 0} \int \rho_\epsilon(l) dl=
\int \lim_{\epsilon \to 0} \rho_\epsilon(l) dl=1.
$$
This shows that $p$ is a density point of $X$.
\end{pf}

We say that a $k$-parameter $F$ satisfies the Tsujii transversality
at $p_0$ if there exists a line through $p_0$
along which the one-parameter Tsujii transversality
condition is satisfied.  In other words, there exists
an affine map $L:B^1 \to B^k$ such that $L(t_0)=p_0$ for some $t_0 \in
\inter B^1$ and such that the induced
one-parameter family $F^L$ defined by $f^L_t=f_{L(t)}$
is Tsujii transverse at the parameter $t_0$.

By linearity of (\ref {eq}) with respect to $v$,
if $(F,p_0)$ is Tsujii transverse
then all lines passing through $p_0$ are Tsujii transverse except
the lines parallel to a certain codimension-one space of $\R^k$.

\begin{lemma} \label {Tt}

Let $F$ be a $C^2$ $k$-parameter family of unimodal maps.  Assume
$(F,p_0)$ satisfies the Tsujii conditions.
Then $p_0$ is a density point of
parameters $p$ for which $(F,p)$ satisfies the Tsujii conditions and
for which $f_p$ is WR.

\end{lemma}

\begin{pf}

If $F$ is Tsujii transverse at $p_0$ then it
is Tsujii transverse along almost every line through $p_0$.
Along such a line it is a density
point of parameters satisfying the Tsujii conditions and WR.
The result follows from Proposition \ref {dimension}.
\end{pf}

\subsubsection{Tsujii transversality and hybrid lamination}
\label {transv cond}

Let us take a closer look at the Tsujii transversality
for an analytic $F$.  Let $f_p=f$.

Assuming the summability condition,
\begin{equation} \label {summability}
\sum_{k=0}^\infty \frac {1} {|Df^k(f(0))|}<\infty
\end{equation}
(in particular if $f$ is CE), let
$$
\nu_f(v)=\sum_{k=0}^\infty \frac {v(f^k(0))} {Df^k(f(0))}=v(0)+
\sum_{k=1}^\infty \frac {v(f^k(0))} {Df^k(f(0))}
$$
be a functional
defined on continuous vector fields $v$ on the interval.

\begin{lemma}

If $f$ satisfies the summability condition then there exists an even
polynomial vector field $v$, with $v(-1)=v(1)=0$ and
such that $\nu_f(v) \neq 0$.

\end{lemma}

\begin{pf}

Let $S=\sum |Df^k(f(0))|^{-1}$.  Let $\epsilon$ be so small that
$$
\sum_{\ntop
{k>0}
{f^k(0) \in (-\epsilon,\epsilon)}}
\frac {1} {|Df^k(f(0))|}<1/3.
$$

Let $v$ be an even polynomial vector field satisfying $v(-1)=v(1)=0$,
\begin{align*}
|v(x)|<2, \quad &\text {for} \quad x \in I,\\
v(x)>1, \quad &\text {for} \quad x \in (-\epsilon/2,\epsilon/2),\\
|v(x)|<\frac {1} {10S}, \quad &\text {for} \quad x \in
I \setminus (-\epsilon,\epsilon).
\end{align*}
Then $\nu_f(v)>1-2/3-1/10>0$.
\end{pf}

\begin{lemma}

The kernel of $\nu_f$ intersected with $T \AAA^\R_a$ is the tangent space to
the hybrid class of $f$.

\end{lemma}

\begin{pf}

By the previous lemma, $\nu_f$ is non-trivial over $T \AAA^\R_a$,
so the above intersection
is a closed codimension-one subspace of $T \AAA^\R_a$.
So it is enough to show that if $v$ is tangent then
$\nu_f(v)=0$.  Assuming that $v$ is tangent, consider an analytic
family $f_t$ contained in the hybrid class of $f$,
such that $f_0=f$ and
$$
\left . \frac {d} {dt} f_t \right |_{t=0}=v.
$$
It is remarked in \cite {ALM} that
$$
\alpha_{n+1}=Df^n(f(0)) \sum_{k=0}^n \frac {v(f^k(0))}
{Df^k(f(0))}=Df^n(f(0)) \nu_f(v)
$$
is precisely
$$
\left . \frac {d} {dt} f_t^{n+1}(0) \right |_{t=0}.
$$
Moreover, $t \mapsto f_t^{n+1}(0)$ are holomorphic functions of the
complex parameter $t$, taking values in $\Omega_a$, and whose domain is
some definite neighborhood of $0$.  It follows by Cauchy estimates on the
derivative that this sequence is bounded
independently of $n$.  By the summability condition
(\ref {summability}), $|Df^n(f(0))| \to \infty$, so we
have necessarily $\nu_f(v)=0$.
\end{pf}

\begin{rem}

It is shown in \cite {ALM} that the sequence $\alpha_n$ is not only
bounded (for tangent vector fields $v$),
but that the vector field defined on the orbit of the
critical value by $w(f^k(0))=\alpha_k$, $k>0$,
extends to a quasiconformal vector field on $\C$.

\end{rem}

So Tsujii transversality can be interpreted for such a map
(satisfying the summability condition (\ref {summability}))
as transversality of the family to the hybrid class of $f_p$.

Since for maps with negative Schwarzian derivative CE
implies the BCE and that all periodic orbits are repelling,
we can conclude from Theorem A, Lemma \ref {analytic}
and this discussion the following result:

\begin{lemma} \label {standard conditions}

If $F$ is a non-trivial $k$-parameter analytic
family of $S$-unimodal maps then almost every parameter
is regular or satisfies the Tsujii conditions.

\end{lemma}

\subsection{Estimates of density in perturbed families}

Let $K$ be the space of
$C^2$ $k$-parameter families of unimodal maps (without, naturally, the
hypothesis of negative Schwarzian derivative).

Let $X \subset K \times B^k$ be the set of $(F,p)$ such that either
$f_p$ is regular or satisfies the Tsujii conditions and WR.
For $F \in K$, let $X_F=\{p \in B^k|(F,p) \in X\}$.

Let $Y \subset B^k$ be measurable with $|Y|>0$.  We define the density of
$X$ along $F$ on $Y$ as
$$
d(F,Y)=\frac {|Y \cap X_F|} {|Y|}.
$$

Instead of defining the classical infinitesimal density:
$$
\liminf_{\epsilon \to 0} d(F,B_\epsilon(p))
$$
we will need to consider the stability of the density with respect to
perturbations of $F$.  With this in mind we introduce two parameters.
Let
$$
D^-(F,p)=\liminf_{\ntop
{\hat F \to F} {\hat f_p=f_p}}
\liminf_{\epsilon \to 0} d(\hat F,B_\epsilon(p)),
$$
$$
D^+(F,p)=\liminf_{\epsilon \to 0} \liminf_{\hat F \to F}
d(\hat F,B_\epsilon(p)).
$$

\begin{rem}

Notice that in the definition of $D^-(F,p)$
we only consider families through a fixed map,
while in the definition of $D^+(F,p)$ we do not make
this restriction.

\end{rem}

Theorem A and Tsujii's result give a direct way to estimate $D^-$:

\begin{lemma}

Let $F$ be a non-trivial analytic family of $S$-unimodal maps.
Then for almost every $p \in B^k$, $D^-(F,p)=1$.

\end{lemma}

\begin{pf}

Indeed, by Lemma \ref {standard conditions}, almost every parameter is either
regular or satisfies the Tsujii conditions.  Since the set of
regular maps is $C^2$ open, $D^-(F,p)=1$ at any regular parameter $p$.

Let us show that this still holds for parameters $p$ satisfying the Tsujii
conditions.  Since Tsujii transversality {\it through a fixed CE map}
is clearly an open condition,
if $\hat F$ is any $C^2$ family near $F$ with $\hat f_p=f_p$
then $(\hat F,p)$ also satisfies the Tsujii conditions.
By Lemma \ref {Tt},
$$
\lim_{\epsilon \to 0} d(\hat F,B_\epsilon(p))=1.
$$
Thus $D^-(F,p)=1$.
\end{pf}

However, for measure estimates in perturbed families,
$D^+(F,p)$ is more relevant.  We proceed to discuss
the effect of the interchange
of limits in the definitions of $D^-(F,p)$ and $D^+(F,p)$.

\begin{lemma} \label {well approximable}

In this setting,
$$
D^+(F,p) \geq D^-(F,p).
$$

\end{lemma}

\comm
{
\begin{pf}

Take sequences $F^{(n)}$ and $\epsilon_j$ satisfying
$$
\lim_{j \to \infty} \liminf_{n \to \infty}
d(F^{(n)},B_{\epsilon_j}(p_0))=D^+(F^{(0)},p_0).
$$
and fix a neighborhood $\VV$ of $F^0$.
Let $k_n$, $l_n$ be the sequences from the approximation condition.  Let
$m_n>l_n$ be a sequence such that
$$
\lim_{n \to \infty} d(F^{(m_n)},B_{\epsilon_{k_n}}(p_0))=
D^+(F^{(0)},p_0).
$$
and let $\tilde F$ be as in the approximation condition.
Then it is clear that
$$
\lim_{n \to \infty} d(\tilde F,B_{\epsilon_{k_n}}(p_0))=
D^+(F^{(0)},p_0).
$$
Since we can get $\tilde F$ to approximate $F^{(0)}$ by shrinking
$\VV$, we get $D^+ \geq D^-$.
\end{pf}
}

\begin{pf}

The idea is to construct, arbitrarily near $F$, a family $\tilde F$ with
$\tilde f_p=f_p$ and
$$
\lim_{j \to \infty} d(\tilde F,B_{\epsilon_j}(p))=D^+(F,p),
$$
for some sequence $\epsilon_j \to 0$, which implies $D^+(F,p) \geq
D^-(F,p)$.
To construct $\tilde  F$, we will interpolate $F$ with a certain sequence
$F^{(n)}$ which realizes the limit in the definition of $D^+(F,p)$.

Let $\epsilon_j \to 0$ be a sequence such that
$$
\lim_{j \to \infty} \liminf_{\hat F \to F}
d(\hat F,B_{\epsilon_j}(p))=D^+(F,p).
$$

Passing to a subsequence, we may assume that
$$
\lim_{j \to \infty} \frac {\epsilon_{j+1}} {\epsilon_j}=0.
$$
Let $K_j \subset B_{\epsilon_j}(p) \setminus
\overline {B_{\epsilon_{j+1}}(p)}$ be compact sets such that
\begin{equation} \label {a}
\lim_{j \to \infty} \frac {|\inter K_j|} {|B_{\epsilon_j}(p)|}=1.
\end{equation}

Let $\phi_j:\R^k \to \R$ be a $C^\infty$ function supported in
$B_{\epsilon_j}(p) \setminus \overline {B_{\epsilon_{j+1}}(p)}$ such
that $\phi_j|K_j=1$.

For a sequence $F^{(n)} \to F$, let us define $\tilde F:B^k \times I
\to I$ by
$$
\tilde f_q=f_q+\sum_{j=1}^\infty \phi_j(q) (f^{(j)}_q-f_q).
$$
It is easy to see that for every $\delta>0$ there exists
a sequence $\delta_n>0$, $n \geq 1$, such that, if
$\|F^{(n)}-F\|_{C^2}<\delta_n$ then
$\|\tilde F-F\|_{C^2}<\delta$ (and in particular $\tilde F$ is $C^2$).
In other words, if $F^{(n)} \to F$ sufficiently fast then
$\tilde F$ is $C^2$ and close to $F$ in the $C^2$ topology.

Notice that $\tilde F$ interpolates $F$ and the sequence $F^{(n)}$
in such a way that inside each
$B_{\epsilon_n}(p)$, $\tilde f_p=f^{(n)}_p$ for $p$ in $\inter K_n$.  Thus,
\begin{equation} \label {a2}
X_{\tilde F} \cap K_n=X_{F^{(n)}} \cap K_n.
\end{equation}

Fix $\delta>0$ and select $F^{(n)}$ such that
\begin{equation} \label {a1}
\lim_{n \to \infty} d(F^{(n)},B_{\epsilon_n}(p))=D^+(F^{(0)},p)
\end{equation}
and moreover $\|F^{(n)}-F\|_{C^2}<\delta_n$, so that $\|\tilde
F-F\|_{C^2}<\delta$.  By (\ref {a}), (\ref {a2}), and (\ref {a1}),
$$
\liminf_{\epsilon \to 0} d(\tilde F,B_\epsilon(p)) \leq \lim_{n \to \infty}
d(\tilde F,B_{\epsilon_n})=\lim_{n \to \infty}
d(F^{(n)},B_{\epsilon_n})=D^+(F,p).
$$
Making $\delta \to 0$, $\tilde F$ converges to $F$ and
we obtain $D^+(F,p) \geq D^-(F,p)$.
\end{pf}

\subsection{Proof of Theorem B}

Let $F$ be a non-trivial analytic family of $S$-unimodal maps.
Then almost every parameter
satisfies $D^-(F,p)=1$.  Hence, for almost every $p$ we have
$D^+(F,p)=1$.

Fix $\epsilon>0$.  Let $p \in B^k$ be such that $D^+(F,p)=1$.
By definition of $D^+$ there exists a sequence of balls $U^n(p)$
centered at $p$ and converging to $p$, and neighborhoods $\VV^n(p) \subset K$
of $F$ such that if $\tilde F \in \VV^n(p)$ then
$$
d(\tilde F,U^n(p)) > 1-\epsilon/2.
$$

By Vitali's Lemma, there exist sequences $p_j$, $n_j$ such that
$U^{n_j}(p_j)$ are disjoint and $|\cup U^{n_j}(p_j)|=1$.  Let $m$ be
such that $\cup_{j=1}^m U^{n_j}(p_j)>1-\epsilon/2$.  Let $\VV=\cap_{j=1}^m
\VV^{n_j}(p_j)$.  Then if $\tilde F \in \VV$, $d(\tilde F,
B^k) \geq 1-\epsilon$.
If $F^{(n)} \to F$ in the $C^2$ topology then $F^{(n)} \in \VV$ for $n$
large enough and the set of parameters for $F^{(n)}$
which are either regular or satisfy the Tsujii conditions and Weak
Regularity have measure at least $1-\epsilon$, as required.

Moreover, considering the sequence $F^{(n)} \equiv F$, we conclude that
almost every parameter for $F$ is Weakly Regular, hence
the last claim of Theorem B.

\subsection{Proof of Theorem C (Smooth Dichotomy)}

By Proposition \ref {strongdense} non-trivial analytic
families are dense among
$C^n$ families of $S$-unimodal maps, $n=3,...,\infty$.
Theorem B implies that for all $\epsilon$ the set $D_\epsilon$
of $C^n$ families of
$S$-unimodal maps for which the set of bad parameters (not regular or
BCE, CE, SE and WR) has measure less then $\epsilon$, contains a
neighborhood of all non-trivial analytic families, that is, an
open and dense set.  Therefore $\cap D_{1/2^n}$ is a residual set.
Clearly any family in $\cap D_{1/2^n}$ satisfies the stated dichotomy.

\appendix

\section{Quasisymmetric robustness of Collet-Eckmann and polynomial
recurrence}

The aim of this Appendix is to sketch a proof of Theorem \ref
{capacity}.
This proof is similar in strategy to the one of the main results of
\cite {AM}, however non-trivial modifications are needed.  To avoid too
much intersection, this will be a concise exposition
concentrated mainly on the new steps needed for this improvement:
the reader can find a full proof of this result
in \cite {AM3}.

\subsection{Quasisymmetric maps}

\subsubsection{Quasisymmetric reparametrization}

Let now $H$ be an arbitrary but fixed
$\hat \gamma$-quasisymmetric map from $B^1$ to 
the parameter space of the quadratic family.  To prove Theorem \ref
{capacity}, it will be enough to show that almost every $t \in B^1$
correspond under $H$ to a parameter of the quadratic family
which is either regular or satisfies
the Collet-Eckmann and Polynomial Recurrence conditions.

{\it From now on, all mentions to parameter space
will (unless explicitly stated otherwise) refer to the above
reparametrization.}

\subsubsection{Quasisymmetric capacities}

The $\g$-capacity of a set $X \subset \R$ in an interval $I$ is defined as
follows:
$$
p_\g(X|I)=\sup \frac {|h(X \cap I)|} {|h(I)|}
$$
where the supremum is taken over all $\g$-qs maps $h:\R \to \R$.

Notice that if $I^j$ are
disjoint subintervals of $I$ and $X \subset \cup I^j$ then
$$
p_\g(X|I) \leq p_\g(\cup_j I^j|I) \sup_j p_\g(X|I^j).
$$

\subsection{Sequence of first return maps}

The statistical analysis of \cite {AM} concerns mainly the
following objects: we are given a unimodal map (which we will assume
finitely renormalizable and with a recurrent critical point)
$f:I \to I$ and a sequence
of nested intervals $I_n \subset I$.  The inductive
relation between the $I_n$ is
as follows: the domain of the first return map
$R_n$ to $I_n$ consists of countably
many intervals $\{I^j_n\}_{j \in \Z}$,
with the convention that $0 \in I^0_n$ (the central component), and we let
$I^0_n=I_{n+1}$.

The special sequence of intervals $I_n$ that we consider is called the
principal nest, see \cite {puzzle}.  Since we assume $f$ to be finitely
renormalizable, there exists a smallest symmetric
interval $T \subset I$ which is
periodic (say, of period $m$).  For the principal nest,
$I_1=[-p,p]$, where $p$
is the orientation reversing fixed point of $f^m:T \to T$.
A level $n$ of the principal nest is
called central if $R_n(0) \in I_{n+1}$.
Let us say that $f$ is a {\bf simple map} if its
principal nest has at most finitely many central levels.

Each non-central branch of $R_n$ is a diffeomorphism onto $I_n$.  Let us
introduce some convenient notation related to the
iteration of the non-central branches of $R_n$.
Let $\Omega$ be the set of finite sequences of non-zero integers (the empty
sequence is included),
an element of $\Omega$ is denoted $\d=(j_1,...,j_m)$.  If $\d \in \Omega$
has length $|\d|=m$, we denote $R_n^\d$
the branch of $R_n^{|\d|}$ with combinatorics $\d$, that is, the domain of
$R_n^\d$ is the set
$$
I^\d_n=\{x \in I|R_n^{k-1}(x) \in I^{j_k}_n,1 \leq k \leq m\}.
$$
We let $C^\d_n=(R_n^\d)^{-1}(I_{n+1})$.

Let us denote by $L_n$ the first landing map from $I_n$ to $I_{n+1}$.
This map relates easily to $R_n$ using
the above description: the domain of $L_n$ is $\cup C^\d_n$, and
$L_n|C^\d_n=R^\d_n$.  The reader should think of $L_n$ as a high iterate of
$R_n$.
This leads to the following inductive relation between return maps:
$R_{n+1}=L_n \circ R_n|I_{n+1}$.

The return time of a point $x$ belonging to an interval
$I^j_n$ is denoted by $r_n(x)$ (or $r_n(j)$, since it does not depend on $x
\in I^j_n$), that is,
$R_n|I^j_n=f^{r_n(j)}$.  The landing time is denoted by $l_n(x) \equiv
l_n(j)$.  The combinatorics at level $n$
of a point $x$ is denoted $\d^{(n)}(x)$, so that $x \in
C^{\d^{(n)}(x)}_n$.  Let $j^{(n)}(x)$ be such that $x \in
I^{j^{(n)}(x)}_n$.  We let $\tau_n=j^{(n)}(R_n(0))$, so that $R_n(0)
\in I^{\tau_n}_n$.  The return time of the critical point is denoted
$v_n=r_n(0)$.  Let $s_n=|\d^{(n)}(R_n(0))|$.

Notice that
$I_{n+1}=R_{n-1}^{-1}(C^\d_{n-1})$ for some $\d$.
The interval $\tilde I_{n+1}=R_{n-1}^{-1}(I^\d_{n-1}) \subset I_n$
is a big neighborhood of $I_{n+1}$ which will be useful later.
This choice of neighborhood is particularly good for simple maps,
and it turns out that in this case $\tilde I_{n+1}$ is still much smaller
than $I_n$ for big $n$.

\subsubsection{Phase-parameter relation}

The starting point of \cite {AM} are two theorems of Lyubich describing
the (unreparametrized)
parameter space of the quadratic family:
infinitely renormalizable maps have zero Lebesgue measure \cite {regular}
and almost every finitely renormalizable non-regular map is simple
\cite {parapuzzle}.
We will need the following remark of \cite {ALM}:
Lyubich's proof actually allows one to conclude that the
set of regular or simple maps has full measure after any quasisymmetric
reparametrization.

In view of those results, Theorem \ref {capacity} is reduced to proving that
the set of parameters which are Collet-Eckmann and polynomially recurrent
have full measure ({\it after reparametrization by $H$}) among simple maps.
From now on we exclude non-simple maps
from measure-theoretic considerations, and we will use ``with total
probability'' to refer to a set of parameters with full measure ({\it after
reparametrization by $H$}) among simple maps.

To estimate the probability in the parameter corresponding to a certain
behavior of the $n$-th stage of the principal nest, we make use of the
Phase-Parameter Lemmas of \cite {AM}.
They describe how the partition of the phase space induced by
return and landing maps $R_n$ and $L_n$ induce parameter partitions of
certain parameter windows $J_n$.

The topological part of the phase-parameter relation is described in the
following:

\begin{thm} \label {top phpa}

For each non-renormalizable quadratic map $f$ with a recurrent critical
point, there exists a sequence of parameter intervals $\{J_n\}$ such that:

\begin{enumerate}

\item $J_n$ is the maximal interval containing $f$ such that for all $g \in
J_n$, there exists a continuation $I_{n+1}[g]$ of $I_{n+1}$ with the
``same combinatorics'' in the following sense.
There exists a continuous family of homeomorphisms
$h_n[g]:I \to I$, $g \in J_n$ which is equivariant with respect to the
actions of $g|(I \setminus I_{n+1}[g])$ and $f|(I \setminus I_{n+1})$, so
that if $x \in I \setminus I_{n+1}[f]$ then $g \circ h_n[g](x)=h_n[g] \circ
f(x)$.

\item There exists a homeomorphism $\Xi_n:I_n \to
J_n$ such that $\Xi_n(C^\d_n)$ is the set of all $g \in J_n$ such that
$R_n[g](0) \in h_n[g](C^\d_n)$.

\end{enumerate}

\end{thm}

This result follows immediately from the Topological Phase-Parameter
relation for the unreparametrized quadratic family (Theorem 2.2
of \cite {AM}), since the reparametrization
is a homeomorphism.

In words, the sequence $J_n$ in Theorem \ref {top phpa}
denotes the maximal interval containing $f$
where we can consider a continuation of $I_n$ (recall that the
boundary of $I_n$ is preperiodic), and such that the first return map
to this continuation does not change combinatorics, so that its domain
changes continuously.  When we change the
map $g$ inside the interval $J_n$, the critical value of $R_n[g]$ varies
inside the interval $I_n[g]$ ``properly'', that is, moves from one boundary
point to the other.  In doing so, it goes through the partition induced by
the $C^\d_n$ in a well behaved (``monotonic'') way: it goes through
each member of the partition exactly once, and thus defines a partition in
the parameter interval $J_n$,
corresponding topologically to the partition in the phase interval $I_n$.
Theorem \ref {top phpa}
thus establishes that the ``diagonal'' motion of the critical
value and the ``horizontal'' motion of the partition of the phase space are
``transversal''.  This is indeed how the proof of Lyubich goes (using
complex analysis).  This result can also be established using the
Milnor-Thurston's combinatorial theory of unimodal
maps together with the monotonicity property of the quadratic family.

The next component of the phase-parameter relation is a quantitative
estimate on the regularity of the phase-parameter homeomorphisms $\Xi_n$. 
While the topological part is based on a very general transversality
argument, the quantitative part depends on the delicate geometric estimates
of Lyubich.

We let $J^{\tau_n}_n=\Xi_n(I^{\tau_n}_n)$.  The correspondence
$\Xi_n$ is uniquely defined if restricted to $K_n=I_n \setminus \cup
C^\d_n$.  More importantly, it is quasisymmetric if restricted
to certain subsets of $K_n$.  To make this precise, let
$K^\tau_n=K_n \cap I^{\tau_n}_n$ (forgetting information outside
$I^{\tau_n}_n$) and
$\tilde K_n=I_n \setminus (\cup I^j_n \cup \tilde I_{n+1})$ (forgetting
information inside each $I^j_n$ and also inside $\tilde I_{n+1}$).

\begin{thm}

Let $f$ be a simple map.  Then, for all
$\g=(1+\delta) \hat \gamma>\hat \gamma$, there exist
$n_0>0$ such that for all $n>n_0$,

\begin{description}

\item[PhPa1] $\Xi_n|K^\tau_n$ is $\g$-qs;

\item[PhPa2] $\Xi_n|\tilde K_n$ is $\g$-qs;

\item[PhPh1] $h_n[g]|K_n$ is $1+\delta$-qs for all
$g \in J^{\tau_n}_n$;

\item[PhPh2] $h_n[g]|\tilde K_n$ is
$1+\delta$-qs for all $g \in J_n$.

\end{description}

\end{thm}

This theorem is a straightforward consequence of the Phase-Parameter
relation for the unreparametrized quadratic family
(Theorem 2.3 of \cite {AM}).
While in \cite {AM} the quasisymmetric constants in PhPa1 and PhPa2
could be taken arbitrarily close to
$1$ (the unreparametrized case corresponds to taking
$H=\id$, that is, $\hat \gamma=1$)
for deeper levels of the principal nest, this does not hold here due
to introduction of reparametrization, which multiply all phase-parameter
constants by $\hat \gamma$ (notice that PhPh1 and PhPh2 are estimates which
do not depend on reparametrization, so we can still choose constants close
to $1$).  This will be the source of
many difficulties addressed in this Appendix.

\subsection{The statistical argument}

For the remaining of this Appendix we fix some constant
$\g>\hat \gamma$, and we will start our consideration with levels
of the principal nest where the
reparametrized phase-parameter relation is already $\g$-qs.
We will also need some very large constants $\b<b$ which depend only
on $\g$ (the relation can be computed explicitly following the proof,
in particular, $\b$ should be at least so big that $\b^{-1}$
is a lower bound on the H\"older constant of $\g$-qs maps).
We let $a=b^{-1}$ and $\a=\b^{-1}$.

From now on we will always estimate the $\g$-capacity of bad sets in the
phase space.  To conclude results for the parameter we will
use the following variation of the Borel-Cantelli Lemma (this is
Lemma 3.1 of \cite {AM}).

\begin{lemma} \label {measure argument}

Let $X \subset \R$ be a measurable set such that for each
$x \in X$ there is a sequence
$D_n(x)$ of nested intervals converging to $x$
such that for all $x_1,x_2 \in X$
and any $n$, $D_n(x_1)$ is either equal or disjoint to $D_n(x_2)$.  Let
$Q_n$ be measurable subsets of
$\R$ and $q_n(x)=|Q_n \cap D_n(x)|/|D_n(x)|$.  Let $Y$ be
the set of $x$ in $X$ which belong to finitely many $Q_n$.
If $\sum q_n(x)$ is finite for almost any $x \in X$ then $|Y|=|X|$.

\end{lemma}

In practice, the $D_n$ will be the parameter windows defined before (either
$J_n$ or $J^{\tau_n}_n$), and
$Q_n$ will be certain subsets of $J_n$ or $J^{\tau_n}_n$
corresponding (under the phase-parameter map) to branches
of the return map (in the case of $J_n$) or landings (in the case of
$J^{\tau_n}_n$),
whose behavior we want to avoid.  We will then show that such bad events
have summable $\g$-capacity in the phase space, which will yield the
conclusion for Lebesgue measure of the parameter using PhPa1 (for landings)
or PhPa2 (for returns).

\subsubsection{A simple application: torrential decay of geometry}

We will now illustrate the use of Lemma \ref {measure argument} and the
phase-parameter relation with an estimate on the decay of geometry.
More precisely, we will consider the {\bf scaling factor}
$$
c_n=\frac {|I_{n+1}|} {|I_n|}.
$$
The scaling factor is a particularly important parameter in the subsequent
analysis: all statistical estimates that follow will be related to $c_n$.

One initial information on the scaling factors is provided by the following
result of Lyubich:

\begin{thm}[see \cite {attractors}] \label {attractorstheorem}

If $f$ is simple than there exists $C>0$, $\lambda<1$ such that $c_n<C
\lambda^n$.

\end{thm}

We will now show that, with total probability, the decay of $c_n$ is much
faster than exponential.  To express this decay, let
us consider the tower function defined by recursion $T(1)=2$,
$T(n+1)=2^{T(n)}$.
We will show that, with total probability, the $c_n$
decrease torrentially to $0$, that is, there exists $k>0$ such that
$c_n^{-1}>T(n-k)$ for $n$ big enough.  More precisely, we will show that
$c_{n+1}^{-1}$ behaves as an exponential of (a bounded power of) $c_n^{-1}$.

This very fast decay
implies that the landing map to $I_{n+1}$ is
essentially a very high iterate
of the return map to $I_n$ (since it takes a long time to hit a very
small interval). This very high iteration time will allow us to
conclude that the characteristics (say, return time) of each level
tend to be better behaved than in the
previous one due to fast convergence to some average (some kind of Law of
Large Numbers).  The fact that we must deal with qs-capacity instead of
Lebesgue measure will essentially reflect in the presence of
errors terms (whose size depend on $\hat \gamma$) in certain exponents
in the above description.

In order to estimate $c_n$, we first consider the related quantity
$s_n=|\d^{(n)}(R_n(0))|$, which denotes the number of times the critical
orbit visits $I_n$ before hitting $I_{n+1}$.

If the critical orbit behaved as a
sequence of random points (uniformly distributed with respect to Lebesgue),
the expectation of this first hitting time should be $c_n^{-1}$.  More
relevant for us, the distribution of the first hitting time (for the random
model) should be concentrated about $c_n^{-1}$: with large probability (say,
less than $2^{-n}$), the
first hitting time is in some ``neighborhood'' of $c_n^{-1}$ (say, $[4^{-n}
c_n^{-1}, 4^n c_n^{-1}]$).  The
corresponding statement for our actual dynamical system
is that the distribution of $|\d^{(n)}(x)|$,
with respect to Lebesgue measure on $x \in I_n$
is concentrated around $c_n^{-1}$, which can be easily checked
by the reader: the estimates are not significantly affected in the
non-random case.

However, due to the nature of the phase-parameter relation, we must
estimate the distribution of $|\d^{(n)}(x)|$ in terms of capacities.  This
will affect drastically the estimates.  To understand why, keep in mind that
$\g$-qs maps are only H\"older (with some constant bounded from below by
$\b^{-1}$), so they can potentially distort the
logarithm of the ratio between $I_{n+1}$ and $I_n$ by such a constant.
Aside from this problem, the information we need can be computed quite
easily and is summarized below.

\begin{lemma} \label {estimate on m}
  
With total probability, for all $n$ sufficiently big we have
\begin{align}
\tag{1}&
p_{2 \g}(|\d^{(n)}(x)| \leq k|I_n)<k c_n^\a,\\
\tag{2}&
p_{2 \g}(|\d^{(n)}(x)| \geq k|I_n)<e^{-k c_n^\b}.
\end{align}
We also have
\begin{align}
\tag{3}&
p_{2 \g}(|\d^{(n)}(x)| \leq k|I^{\tau_n}_n)<k c_n^\a,\\
\tag{4}&
p_{2 \g}(|\d^{(n)}(x)| \geq k|I^{\tau_n}_n)<e^{-k c_n^\b}.
\end{align}

\end{lemma}

This lemma corresponds to Lemma 4.2 of \cite {AM}.

The phase-parameter lemmas
(specially PhPa1) allow us to transfer the last pair of estimates to the
parameter space: for $n$ sufficiently big,
(Lebesgue) most parameters in $J^{\tau_n}_n$ satisfy
$$
c_n^{-\a/2}<s_n<c_n^{-2\b}.
$$
Here `most' means that the
complement has probability bounded by $c_n^{\a/3}$.  But $c_n$ (and thus
$c_n^{\a/3}$) decays exponentially for every simple map (by Theorem \ref
{attractorstheorem}).  So $\sum c_n^{\a/3}<\infty$ and we
are able to apply Lemma
\ref {measure argument} to obtain the following:

\begin{lemma} \label {growth of s_n}

With total probability, for $n$ sufficiently big we have
$$
c_n^{-\a/2}<s_n<c_n^{-2\b}.
$$
\end{lemma}

This lemma corresponds to Lemma 4.3 of \cite {AM}.

\begin{rem} \label {c_n t}
This result implies easily torrential decay of $c_n$: $\ln
c_{n+1}^{-1}$ can be easily bounded from below by $K s_n$ for some universal
$K>0$, and thus for big $n$,
$$
c_{n+1}^{-1} \geq e^{c_n^{-\a/3}}.
$$
\end{rem}

\subsection{Derivatives}

We proceed to estimate derivatives of branches of the return map. 
All lemmas in this section can be proved using the
same argument as in \cite {AM}.

The first step is to exclude the possibility of a `too recurrent' 
or `too low' return.  It is analogous to Lemma 4.8 of \cite {AM}, being a
simple application of PhPa2.

\begin{lemma} \label {away from the boundary}

With total probability, the distance between $R_n(0)$ and
$\partial I_n \cup \{0\}$ is at least $|I_n|n^{-\b}$.
In particular $R_n(0) \notin \tilde I_{n+1}$ for all
$n$ large enough.

\end{lemma}

Recall that the {\bf distortion} of a diffeomorphism $\phi$
on an interval $T$ is defined by
$$
\Dist(\phi|T)=\frac {\sup_T |D\phi|}
{\inf_T |D\phi|}.
$$

Lemma \ref {away from the boundary} allows us to start
estimating the distortion of iterates of $f$.
The following estimate corresponds to Lemma
4.9 of \cite {AM}.  It is based on the fact
that the distortion of branches of return
maps is due to the position of the branch with respect to the critical
point.  Using PhPa1, we are able to give polynomial lower
bounds on the distance between the
critical point with respect to non-central branches, which are valid with
total probability.

\begin{lemma} \label {dist}

With total probability, for $n$ big enough and $j \neq 0$
$$
\Dist(f|I^j_n) \leq n^\b.
$$

\end{lemma}

The following estimate is analogous to Lemma 4.10 of \cite {AM}.  It is
based on the previous one and the observation that
return branches are torrentially expansive in average (from the decay 
of geometry).

\begin{lemma} \label {distortion}

With total probability, for $n$ big enough and for all $\d \in \Omega$
$$
\Dist(R^\d_n) \leq n^\b.
$$
In particular, for $n$ big enough,
$|DR_n(x)|>2$ if $x \in \cup_{j \neq 0} I^j_n$.

\end{lemma}

Lemma \ref {distortion} gives estimates of derivatives
under iterates of $R_n$.  To obtain estimates of derivatives under iterates
of $f$, we will need the following very general result of
Guckenheimer which shows that quadratic
maps are hyperbolic away from critical points and parabolic points
(this actually generalizes to very general one-dimensional
systems by a result of Ma\~n\'e), see \cite {MvS}.  We state just a
consequence adapted to our particular setting.

\begin{thm} \label {hypthm}

Let $f$ be a quadratic map without non-repelling periodic orbits (in
particular if $f$ is a simple map).
For every $\epsilon>0$, there exists $C>0$, $\lambda>1$
such that if $|f^k(x)|>\epsilon$ for $0 \leq k \leq m$ then
$Df^{m+1}(x)>C \lambda^m$.

\end{thm}

With this information we are now able to give a lower bound on the
derivative of iterates of $f$.  The next lemma is
identical to Lemma 4.11 of \cite {AM}, and is based on the idea that full
returns to sufficiently deep levels cause expansion (from the previous
lemma), while the dynamics
outside a definite neighborhood of the critical point is hyperbolic
(by Theorem \ref {hypthm}).

\begin{lemma} \label {lower bound}

With total probability, if $n$ is sufficiently big and if
$x \in I^j_n$, $j \neq 0$, and $R_n|I^j_n=f^r$, then
for $1 \leq k \leq r$, $|Df^k(x)|>|x| c_{n-1}^3$.

\end{lemma}

\subsection{How to deal with hyperbolicity}

Keeping in mind that our analysis of the statistical properties of the
dynamics of $f$ is made in terms of the induced return maps $R_n$, we see
that in order to estimate the hyperbolicity along the critical orbit (to
obtain the Collet-Eckmann condition) we must have a convenient way
to quantify the
hyperbolicity of (for instance)
non-central return branches.
To do so, for $j \neq 0$, we define the quantity
$$
\lambda_n(j)=\inf_{x \in I^j_n} \frac
{\ln |DR_n(x)|} {r_n(j)}.
$$
We let $\lambda_n=\inf_{j \neq 0} \lambda_n(j)$.

To analyze the behavior of $\lambda_n$, we start with the general
information provided by Theorem \ref {hypthm}.
Coupled with exponential upper bounds on distortion for returns
(which competes with torrential expansion of each non-central branch
from the decay of $c_n$), the hyperbolicity of $f$ in
the complement of $I_{n+1}$ immediately implies the following estimate
(identical to Lemma 7.9 of \cite {AM}).

\begin{lemma}

With total probability, for all $n$ sufficiently big, $\lambda_n>0$.

\end{lemma}

The ``minimum hyperbolicity'' $\liminf \lambda_n$ of the parameters we
will obtain will in fact be positive, as it follows from one of the   
properties of Collet-Eckmann parameters (uniform hyperbolicity on periodic
orbits), together with our estimates on distortion.

Our strategy however is not to show that the minimum hyperbolicity
is positive, but that the typical value of $\lambda_n(j)$
stays big as $n$ grows (and is in fact bigger than $\lambda_{n_0}/2$
for $n>n_0$ big).  In this sense,
it is convenient to think of $\lambda_n(j)$ as a random
variable whose distribution we are interested in.

There is an inductive
relation between the random variables $\lambda_n(j)$ for different values of
$n$: this is related to the fact that if $R_n(I^j_{n+1})
\subset C^\d_n$, $\d=(j_1,...,j_m)$, we have $R_{n+1}|I^j_{n+1}=L_n|C^\d_n
\circ R_n|I^j_{n+1}$.  The hyperbolicity of the ``landing part''
$L_n|C^\d_n$ is essentially a weighted sum
\begin{equation} \label {weighted sum}
\frac {\sum_{i=1}^m \lambda_n(j_i) r_n(j_i)} {\sum_{i=1}^m r_n(j_i)}.
\end{equation}

So if the ``return part'' $R_n|I^j_{n+1}$ does not carry a big weight on the
computation of $\lambda_{n+1}(j)$ (outside a set of branches with small
$\g$-qs capacity), we can think of $\lambda_{n+1}(j)$ as
distributed according to the weighted sum (\ref {weighted sum}).
This turns out to be the case as the return part does not affect much
the denominator (time) and does not
have a bad effect on the numerator (derivative).
Indeed, in the next section we will see that the return time of
$R_n|I^j_{n+1}$ (given by $v_n$) is much smaller (of order $c_{n-1}^{-1}$)
than the total return time $R_{n+1}|I^j_{n+1}$ (of order $c_n^{-1}$). 
Moreover, if $I^j_{n+1}$ is outside a small neighborhood of $0$,
$|DR_n|I^j_{n+1}|$ is bigger than $1$.

Since we also have to estimate the
hyperbolicity of truncated branches (as the Collet-Eckmann condition is a
condition along the full critical orbit, and not only at full returns),
it will not be enough to just obtain that
the distribution of $\lambda_n(j)$ is
concentrated around some value bigger than $\lambda_{n_0}/2$.
In order to state exactly what kind of hyperbolicity estimate we need,
it is convenient to introduce a certain class of branches: good returns.

We define the set of good returns
$G(n_0,n) \subset \Z \setminus \{0\}$, $n_0,n \in \N$,
$n \geq n_0$ as the set of all $j$ such that

\begin{description}

\item[G1] (hyperbolic return)
$$
\lambda_n(j) \geq \lambda_{n_0} \frac {1+2^{n_0-n}} {2},
$$
  
\item[G2] (hyperbolicity in truncated return)
for $c_{n-1}^{-3/(n-1)} \leq k \leq r_n(j)$ we have
$$
\inf_{I^j_n} \frac {\ln |Df^k|} {k} \geq
\lambda_{n_0} \frac {1+2^{n_0-n+1/2}} {2}-c_{n-1}^{2/(n-1)}.
$$
  
\end{description}

Of course we still have to show that the set of returns which fail to be
good has small $\g$-qs capacity.  In order to do so, we will construct
explicitly a class of branches whose complement has small
$\g$-qs capacity and then show that this class of branches is contained in
good branches (see Lemma \ref {very good
is good}).  Before doing so, we must first estimate
the distribution of return times, since they have an important role in the
computation of $\lambda_n(j)$.

\subsection{Distribution of return and landing times}

To estimate the distribution of return and landing times, it is convenient
to also think of $r_n(j)$ and $l_n(j)$ as 
``random variables'' which are related by some simple rules: if
$\d=(j_1,...,j_m)$ then $l_n(\d)=\sum_{i=1}^m r_n(j_i)$ and
$r_{n+1}(j)=v_n+l_n(\d)$ where $R_n(I^j_{n+1}) \subset C^\d_n$.  In
particular, since the distribution of $|\d^{(n)}|$ is concentrated around
$c_n^{-1}$ which is torrentially big, the random variable
$l_n$ behaves like a very large sum of random variables distributed as
$r_n$.  On the other hand, $r_{n+1}$ should have distribution
approximately like $l_n$ itself, once we show that $v_n$ does not make an
important contribution.

The main tool to do the actual analysis is to prove first a
Large Deviation Estimate for $r_n$ using only the torrential decay of $c_n$,
and then show that such estimate leads to much more precise control of the
subsequent levels.

Since the transition between different levels
introduces some distortion (although torrentially small), we are forced to
deal with a sequence of quasisymmetric constants in our estimates: instead
of just estimating $\g$-qs capacities for some fixed $\g$,
we must consider a sequence $\g_n=\g (n+1)/n$ and
$\tg_n=\g (2n+3)/(2n+1)$.  The basic idea is that
control of the distribution of $r_n$ with respect
to $\g_n$-capacities will provide control of the distribution of
$l_n$ with respect to $\tg_n$
capacities which in turn will allow to estimate the distribution of
$r_{n+1}$ with respect to $\g_{n+1}$ capacities.  Notice that $\inf
\g_n=\inf \tg_n=\g$.  (This ideas are introduced in \S 5 of \cite {AM}.)

Although very technical, this part is very similar to the analysis made on
(the several lemmas of)
\S 6 of \cite {AM} (differing only by change of constants),
so we will only state the final estimate which summarizes
the results of that section and provide a short outline of the argument.

\begin{lemma} \label {large times}

With total probability, for all $n$ sufficiently large we have

\begin{enumerate}

\item $p_{\tg_n}(l_n(x)<c_n^{-s}|I_n)<c_n^{\a^2-s} <
c_n^{a-s}$, with $s>0$,

\item $p_{\tg_n}(l_n(x)<c_n^{-s}|I^{\tau_n}_n)<c_n^{a-s}$, with
$s>0$,  

\item $p_{\tg_n}(l_n(x)>c_n^{-s}|I_n)<e^{-c_n^{b-s}}$, with $s>b$,

\item $p_{\tg_n}(l_n(x)>c_n^{-s}|I^{\tau_n}_n)<e^{-c_n^{b-s}}$,
with $s>b$,

\item $p_{\g_n}(r_n(x)<c_{n-1}^{-s}|I_n)<c_{n-1}^{\a^2-s} <
c_{n-1}^{a-s}$, with $s>0$,

\item $p_{\g_n}(r_n(x)>c_{n-1}^{-s}|I_n)<e^{-c_{n-1}^{\sqrt b-s}} <
e^{-c_{n-1}^{b-s}}$ with $s>b$.

\item $c_{n-1}^{-a}<r_n(\tau_n)<c_{n-1}^{-b}$.

\item $c_{n-1}^{-a}<v_n<c_{n-1}^{-b}$.

\item $c_{n-1}^{-a}<\ln (c_n^{-1})<c_{n-1}^{-b}$.

\end{enumerate}

\end{lemma}

\subsubsection{Outline of the proof of Lemma \ref {large times}}

The estimates from below are relatively easy.
Estimates (1) and (2) follow directly from $l_n(\d) \geq
|\d|$ and Lemma \ref {estimate on m}.  Estimate (5) follows from (1) using
the relation between $r_{n+1}$ and $l_n$.  The estimate from below in (8)
follows from (2) and PhPa1, and the estimate from below in (7)
follows from (5) and PhPa2.  The estimate from below on (9) was computed on
Remark \ref {c_n t}.

The estimates from above are much more delicate.  In what follows we
will ignore the difference between $I_n$ and $I^{\tau_n}_n$, since it is not
substantial for the argument.
The key estimate is (6), which says that the tail $p_{\g_n}
(r_n(x)>k)$ decays exponentially fast (in $k$)
with some specific rate (polynomial in $c_{n-1}$).  On the
other hand, decay with {\it some} rate is easy: $f$ is hyperbolic outside
$I_{n+1}$ (see Theorem \ref {hypthm}), so there exists some (small)
$\alpha_n>0$ with
$p_{\g_n}(r_n(x)>k \alpha_n^{-1})<e^{-k}$ for
$k \geq 1$.  This exponential decay implies that
it is very unlikely that a large sequence
$\d=(j_1,...,j_m)$ will have a landing time
$l_n(\d)=\sum_{i=1}^m r_n(j_i)$ much bigger than $m \alpha_n^{-1}$.

From this relation between $r_n$ and $l_n$, we see that
there exists some $\beta_n$ with
$p_{\tg_n}(l_n(x)>k \beta_n^{-1})<e^{-k}$, and moreover we can estimate
$\beta_n$ in terms of $\alpha_n$ and the size of a typical $\d^{(n)}$ (which
is given by a polynomial on $c_n^{-1}$):
$\beta_n^{-1}$ is bounded by a polynomial (this
polynomial error is related to $\g$) on
$\alpha_n^{-1} c_n^{-1}$.
From the relation between $l_n$ and $r_{n+1}$
we obtain an estimate on $\alpha_{n+1}$ in terms of $v_n$ and $\beta_n$,
which we can rewrite in terms of $v_n$, $c_n$ and $\alpha_n$:
$\alpha_{n+1}^{-1}-v_n$ is bounded by some polynomial on
$\alpha_n^{-1} c_n^{-1}$.

Since $p_{\tg_n} (l_n(x)>\beta_n^{-1} c_n^{-1})$ is summable
(by definition of $\beta_n$), it follows that $v_{n+1}-v_n$ is bounded
by a polynomial on
$\alpha_n^{-1} c_n^{-1}$ with total probability (use PhPa1), in particular,
for $n$ big we can bound $v_{n+1}$ with a polynomial on
$\alpha_n^{-1} c_n^{-1}$.

In particular, if $\alpha_n^{-1}>c_n^{-1}$, $\alpha_{n+1}^{-1}$
is bounded by a {\it polynomial} in $\alpha_n^{-1}$.  Although initially
we did not have any control on the value
of $\alpha_n$, we know that
$c_{n+1}^{-1}$ behaves as an {\it exponential} on $c_n^{-1}$ (torrential
growth), so eventually it
catches up with $\alpha_n^{-1}$: for $n$ big, $c_n^{-1}>\alpha_n^{-1}$.

So for $n$ big
$\alpha_n^{-1}$ can be bounded exclusively by a polynomial on
$c_{n-1}^{-1}$ as stated in (6).  This automatically implies the estimate
from above in (7) using PhPa2.  Since $\beta_n^{-1}$ and $v_{n+1}$ are
bounded by a polynomial on $\alpha_n^{-1} c_n^{-1}$
we obtain (3) and (4) and the estimate from above in (8).

Since $f^{v_n}$ expands $I_{n+1}$ to an interval of size at least
$2^{-n}|I_n|$, and the derivative of $f$ is bounded by $4$, we have
$2^n c_n^{-1}<4^{v_n}$, so the estimate from above on
(9) follows from the estimate
from above in (8).

\subsection{Constructing hyperbolic branches}

In this section we show by an inductive process that the great
majority of branches are reasonably hyperbolic (good branches).
In order to do that, in the
following subsection, we define some classes of branches with
`very good' distribution of times and which are not too close
to the critical point.  The definition of `very good' distribution of
times has an inductive
component: they are composition of many `very good' branches of the previous
level.  The fact that most branches are `very good'
is related to the validity of
some kind of Law of Large Numbers estimate.  The inductive definition will
guarantee that the `very good' distribution of times holds in all scales and
allows us to preserve hyperbolicity from one step to the other: very good
branches are good.

\begin{rem}

The several classes of branches that
we will define do not correspond exactly
to the same classes in \cite {AM}, although classes with the same name  
have essentially the same function in the proof.  There are some non-trivial
steps to make this adaptation work,
since the previous proof uses strongly small quasisymmetric
constants.  This will lead to consideration of extra classes below (bad
returns and fast landings).

\end{rem}

\begin{rem}

This section contains the main modifications with respect to \cite {AM}
(precisely the introduction of bad returns and fast landings).  The role of
those modifications is explained in Remark \ref {modifications}.

\end{rem}

\subsubsection{Standard landings}

Let us define the set of standard landings at time $n$,
$LS(n) \subset \Omega$ as the set
of all $\d=(j_1,...,j_m)$ satisfying the following:

\begin{description}

\item[LS1] ($m$ is not too small or large) $c^{-a/2}_n<m<c^{-2b}_n$,

\item[LS2] (No very large times) $r_n(j_i)<c^{-3b}_{n-1}$ for all $i$.

\item[LS3] (Short times are sparse in large enough initial segments)
For $c^{-2b}_{n-1} \leq k \leq m$
$$
\#\{1 \leq i \leq k,\, r_n(j_i)<c^{-a/2}_{n-1}\} <
(6 \cdot 2^n) c^{a/2}_{n-1} k.
$$
  
\end{description}

We also define the set of fast landings at time $n$,
$LF(n) \subset \Omega$  by the following conditions

\begin{description}

\item[LF1] ($m$ is small) $m<c^{-a/2}_n$.

\item[LS2] (No very large times) $r_n(j_i)<c^{-3b}_{n-1}$ for all $i$.
  
\end{description}


It is easy to convince oneself that most landings are standard.  Indeed,
the distribution of $|\d^{(n)}(x)|$ is concentrated
around $c_n^{-1}$ as requested by LS1.  Moreover, branches with
very large times (larger than
$c_{n-1}^{-3b}$) are so few that even a long sequence $(j_1,...,j_m)$ with
$m<c_{n-1}^{-2b}$ is not likely to contain such an event, as required by
LS2.  Finally, the Law of Large Numbers indicates that
a long sequence $(j_1,...,j_m)$ will seldom contain a
proportion of short times much bigger than their frequency as given by Lemma
\ref {large times}, as required by LS3.

Since fast landings are not standard, they must be few.  However, they
correspond to most of the branches which are not standard.  The reason for
this comes from the requirements of LS1, which imposes two conditions (an
upper and a lower bound on $m$).  The upper bound condition is much more
rarely violated (by one exponential order of magnitude) than the lower
bound (just check Lemma \ref {estimate on m}).
Fast landings essentially capture the violations of the lower
bound (LF1).

The actual estimates for the frequency of standard and fast landings are
provided below.  They can be obtained from the estimates of
distribution of return times (contained in Lemma \ref {large times})
following the general lines of Lemma 7.1 of \cite {AM}.
This step is purely dynamical (no further parameter exclusion is made).

\begin{lemma} \label {standard landing}

With total probability, for all $n$ sufficiently big,

\begin{enumerate}

\item $p_{\tg_n}(\d^{(n)}(x) \notin LS(n)|I_n)<c_n^{a/3}/2$,

\item $p_{\tg_n}(\d^{(n)}(x) \notin
LS(n) \cup LF(n)|I_n)<c_n^{n^2}/2$,

\item $p_{\tg_n}(\d^{(n)}(x) \notin LS(n)|I^{\tau_n}_n)<c_n^{a/3}/2$,

\item $p_{\tg_n}(\d^{(n)}(x) \notin LS(n) \cup LF(n)|
I^{\tau_n}_n)<c_n^{n^2}/2$.

\end{enumerate}
  
\end{lemma}

\comm{
\begin{pf}

The proof is immediate from our time estimates, which can be applied
directly to estimate the losses of LS1, LS2, LF1, or under large
deviation form for LS3, following \S \ref {more tree}.
Let us estimate the complement of the sets which do not satisfy some
properties:

\begin{description}

\item [LS1] $c_{n-1}^{2a/5}$,

\item [LS1+LF1] $e^{c_n^b} \ll c_n^{n^2}$,

\item [LS2] $e^{c_{n-1}^{-2b}} \ll c_n^{n^2}$,

\item [LS3] let $q=(6 \cdot 2^n) c_{n-1}^{a/2}$, we get
$$
\sum_{k>c_{n-1}^{-2b}} \left (\frac {1} {2} \right )^{qk} \ll c_n^{n^2}.
$$

\end{description}

This gives immediately $1$ and $2$.  For $3$ and $4$ the same estimates
hold, using (for LS2) the estimate on $r_n(\tau_n)$.
\end{pf}
}
\subsubsection{Very good returns, bad returns and excellent landings}

Define the set of very good returns,
$VG(n_0,n) \subset \Z \setminus \{0\}$, $n_0 \leq n \in \N$ and
the set of bad returns,
$B(n_0,n) \subset \Z \setminus \{0\}$, $n_0 \leq n \in \N$,
by induction as follows.  We let
$VG(n_0,n_0)=\Z \setminus \{0\}$, $B(n_0,n_0)= \emptyset$ and supposing
$VG(n_0,n)$ and $B(n_0,n)$ defined, define the set of excellent landings
$LE(n_0,n) \subset LS(n)$ satisfying the following extra assumptions.   

\begin{description}

\item[LE1] (Not very good moments are sparse in large enough initial
segments)
For all $c^{-2b}_{n-1}<k \leq m$
$$
\#\{1 \leq i \leq k,\, j_i \notin VG(n_0,n)\} <
(6 \cdot 2^n) c^{a^2}_{n-1} k,
$$
\item[LE2] (Bad moments are sparse in large enough initial segments)
For all $c^{-1/n}_{n}<k \leq m$
$$
\#\{1 \leq i \leq k,\, j_i \notin B(n_0,n) \} <
(6 \cdot 2^n) c^n_{n-1} k,
$$
  
\end{description}

We define $VG(n_0,n+1)$ as the set of $j$ such that
$R_n(I^j_{n+1})=C^{\d}_n$ with $\d \in LE(n_0,n)$ and
the extra condition:

\begin{description}

\item[VG] (distant from $0$) The distance of $I^j_{n+1}$ to $0$ is bigger
than $c_n^{n^2}|I_{n+1}|$.

\end{description}

And we define $B(n_0,n+1)$ as the set of $j \notin VG(n_0,n+1)$ such that
$R_n(I^j_{n+1})=C^\d_n$ with $\d \notin LF(n)$.

Very good returns are designed to carry hyperbolicity from level to level:
since they are composed of many very good returns of the previous level
(LE1), and are not too close to $0$ (VG), they should keep most of the
hyperbolicity of level $n_0$ (given by $\lambda_{n_0}>0$).
For this to work, we must control the
distribution of return times of the previous level
inside a very good branch.  The risky situation is the presence of not very
good branches which have a large return time: those are contained in the
bad branches defined above.  It turns out that they can not spoil
the hyperbolicity because they are too few (LE2).  This basic idea will be
carried out in detail through a series of lemmas.

Very good and bad returns can be estimated in an inductive fashion
analogously to the estimate of Lemmas 7.2 and 7.3 of \cite {AM}: initially
all branches are very good and no branches are bad, and as $n$ grows the Law
of Large Numbers indicates that conditions LE1 and LE2 should be rarely
violated so that very good branches should continue to be
frequent and bad branches rare.
This estimate is again purely dynamical.

\begin{lemma} \label {induction step}
With total probability, for all $n_0$ sufficiently big,
\begin{enumerate}
\item $p_{\g_n}(j^{(n)}(x) \notin VG(n_0,n)|I_n)<c_{n-1}^{a^2},$
\item $p_{\g_n}(j^{(n)}(x) \in B(n_0,n)|I_n)<2 c_{n-1}^{2 n}$,
\item $p_{\tg_n}(\d^{(n)}(x) \notin LE(n_0,n)|I_n)<c_n^{2a/5}$,
\item $p_{\tg_n}(\d^{(n)}(x) \notin LE(n_0,n) \cup LF(n)|I_n) <
c_n^{bn}$,
\item $p_{\tg_n}(\d^{(n)}(x) \notin LE(n_0,n)|I^{\tau_n}_n) <
c_n^{2a/5}$.
\end{enumerate}

\end{lemma}

This translates immediately using PhPa2 to a parameter estimate
analogous to Lemma 7.4 of \cite {AM}:

\begin{lemma} \label {lands very good}

With total probability, for all
$n_0$ big enough, for all $n$ big enough (depending on $n_0$),
$\tau_n \in VG(n_0,n)$.

\end{lemma}

Before going on we will need two simple estimates: one is
for the return time of very good branches
and another is for the return time of
branches which are neither very good or bad.  The first of those
estimates is analogous to Lemma 7.5 of \cite {AM},
and follows directly from the definitions of very good and bad branches.

\begin{lemma} \label {very good return time}

With total probability, for all
$n_0$ big enough and for all $n \geq n_0$,
if $j \in VG(n_0,n+1)$ then
$$
m<r_{n+1}(j)<m c^{-4b}_{n-1},
$$
where, as usual, $m$ is such that $R_n(I^j_{n+1})=C^\d_n$ and
$\d=(j_1,...,j_m)$.

\end{lemma}

\begin{lemma} \label {not very good or bad}

With total probability for all $n_0$ sufficiently big, if $n>n_0$,
if $j \notin VG(n_0,n) \cup B(n_0,n)$ then
$r_n(j)<c_{n-1}^{-a/2} c_{n-2}^{-4b}$.

\end{lemma}

\begin{pf}

Indeed, if $j \notin VG(n_0,n) \cup B(n_0,n)$ then $R_{n-1}(I^j_n) \subset
C^\d_{n-1}$ with $\d \in LF(n-1)$.  By definition of fast landing,
$l_{n-1}(\d)<c_{n-1}^{-a/2} c_{n-2}^{-3b}$, so
$$
r_n(j)=v_{n-1}+l_{n-1}(\d)<c_{n-1}^{-a/2} c_{n-2}^{-3b}+c_{n-2}^{-b}.
$$
\end{pf}

At this stage we have most of the tools to show that
almost every parameter is ``Collet-Eckmann at first returns'', that is,
$|Df^{k_n}(f(0))|$
is exponentially big for the sequence $k_n$ of first landings of $f(0)$ in
$I_n$.  To obtain the full Collet-Eckmann condition (exponential growth
for all $k$), we will need to analyze truncations of branches or
landings, that is, we will consider iterates of the type $f^k|I^j_n$
(or $f^k|C^\d_n$)
for $k$ less then the return time $r_n(j)$ (or $l_n(\d)$).

We now show that very good branches are well behaved when
truncated at a reasonably big time.  Here ``well behaved''
means ``spending most of the time in very good branches of
the previous level''.  So if we
are able to control the hyperbolicity of very good branches in some level
we will have a good possibility of controlling
truncated very good branches in
the next level.  This lemma corresponds to Lemma 7.6 of \cite {AM}, but
the proof must be modified, with the use of bad returns and fast landings.

\begin{lemma} \label {very good partial time}

With total probability, for all
$n_0$ big enough and for all $n \geq n_0$,
the following holds.

Let $j \in VG(n_0,n+1)$, as usual let
$R_n(I^j_{n+1}) \subset C^\d_n$ and
$\d=(j_1,...,j_m)$.  Let $m_k$ be biggest possible with
$$
v_n+\sum_{j=1}^{m_k} r_n(j_i) \leq k
$$
(the amount of full returns to level $n$ before time $k$) and let
$$
\beta_k=\sum_{\ntop {1 \leq i \leq m_k,}
{j_i \in VG(n_0,n)}} r_n(j_i).
$$
(the total time spent in full returns to level $n$
which are very good before time $k$)
Then $1-\beta_k/k<c_{n-1}^{a^2/3}$ if $k>c_n^{-2/n}$.

\end{lemma}

\begin{pf}

Let us estimate first the time $i_k$ which is not spent on non-critical
full returns:
$$
i_k=k-\sum_{j=1}^{m_k} r_n(j_i).
$$
This corresponds exactly to $v_n$ plus some incomplete part of the return
$j_{m_{k+1}}$.  This part can be bounded by $c_{n-1}^{-b}+c_{n-1}^{-3b}$ 
(use the estimate of $v_n$ and LS2 to estimate
the incomplete part).

Using LS2 we conclude now that
$$
m_k>(k-c_{n-1}^{-b})c_{n-1}^{3b}>c_n^{-1/n}
$$
so $m_k$ is not too small.

Let us now estimate the contribution $h_k$ from bad full returns $j_i$.  
The number of such returns
must be less than $c_{n-1}^{n/2} m_k$ by LE2 and the estimate on $m_k$.
By LS2 their total time is at most
$c_{n-1}^{(n/2)-3b} m_k<m_k$.

The non very good full returns on the other hand can be estimated by LE1
(given the estimate on $m_k$), they are at most $c_{n-1}^{a^2} m_k$.
So we can estimate the total time $l_k$ of non very good or bad
full returns (with time less then $c_{n-1}^{-a/2} c_{n-2}^{-4b}$ by Lemma
\ref {not very good or bad}) by
$$
c_{n-1}^{a^2}c_{n-1}^{-a/2} c_{n-2}^{-4b} m_k,
$$
while $\beta_k$ can be estimated from below by
$$
(1-c^{a/4}_{n-1}) c^{-a/2}_{n-1} m_k.
$$      

It is easy to see then that $i_k/\beta_k \ll c^{a/5}_{n-1}$,
$h_k/\beta_k \ll c^{a/5}_{n-1}$.
We also have
$$
\frac {l_k} {\beta_k}<2 c_{n-1}^{a^2/2}.
$$
So $(i_k+h_k+l_k)/\beta_k$ is less then $c^{a^2/3}_{n-1}$.
Since $i_k+h_k+l_k+\beta_k=k$ we have $1-\beta_k/k<(i_k+h_k+l_k)/\beta_k$.
\end{pf}

\begin{rem} \label {modifications}

This lemma illustrates the main reason why the original argument of \cite
{AM} must be changed in order to deal with big quasisymmetric constants. 
Indeed, in \cite {AM}, we do not need to split the branches which are not
very good in bad branches and otherwise (fast).  The reason is that in \cite
{AM} the distribution of $r_n(j)$ is concentrated in a much narrower window
around $c_{n-1}^{-1}$
(say, $(c_{n-1}^{-1+2\epsilon},c_{n-1}^{-1-2\epsilon})$). 
In particular, in a large sequence $(j_1,...,j_k)$ (which should be thought
as an initial segment of an excellent landing),
we can estimate the proportion of the total return time
due to very good branches essentially by considering the proportion of very
good branches in the sequence.

In this Appendix, the distribution of $r_n(j)$ is located in a much larger
window $(c_{n-1}^{-a},c_{n-1}^{-b})$.  The risky situation is to have
a large sequence $(j_1,...,j_k)$ with a large proportion of very
good branches, but whose return time is near the bottom of the window
($c_{n-1}^{-a}$), while the not very good branches in the sequence have all
return time near the top ($c_{n-1}^{-b}$).  In this case, the proportion of
the total time due to very good branches could be very small.

The solution given in this Appendix is based on the idea that the not
very good branches {\it with large time} (bad branches)
are really very few: most of the not very good branches are indeed fast. 
Paying attention to this asymmetry, we can indeed prove that in such a
sequence $(j_1,...,j_k)$, most of the total time is due to very good
branches.

This argument (most branches with atypical time are fast)
is based implicitly in the following asymmetry which appeared already in our
first statistical estimate, Lemma \ref {estimate on m}, when we showed that
the distribution of $|\d^{(n)}(x)|$ is concentrated around $c_n^{-1}$:
there is a big difference (one extra exponential)
in the estimates on the upper tail ($\g$-qs capacity
of $\{|\d^{(n)}(x)|>c_n^{-k\b}\}$)
and the lower tail ($\g$-qs capacity of $\{|\d^{(n)}(x)<c_n^{-k\a}\}$).

(Essentially the same problem, with the same solution, appears
in Lemma \ref {cool hyperbolicity}.)

\end{rem}

Now we conclude that very good (that is, most) branches are good, justifying
our previous hints.

\begin{lemma} \label {very good is good}

With total probability, for $n_0$ big enough and for all $n>n_0$,
$VG(n_0,n) \subset G(n_0,n)$.

\end{lemma}

The proof is the same as for Lemma 7.10 of
\cite {AM}, the two main features of very good
branches exploited here are their good distribution of return times and
the condition VG which allows us to avoid drastic losses of
derivative due to starting very close to the critical point.  The argument
is by induction:
first, all very good branches of level $n_0$ satisfy condition G1 of
a good branch, that is, a full return is very hyperbolic
(this follows from the definition of $\lambda_{n_0}$).
Then, supposing that all very
good branches of level $n$ satisfy G1, we conclude that very good
branches of level $n+1$ have enough hyperbolic branches in its composition
(even if truncated) to satisfy both conditions G1 and G2.

\subsubsection{Cool landings}

As we hinted in the last section,
very good branches play the role of building blocks of
hyperbolicity.  We must now show that the critical point spends most of its
time in very good branches.  To do so, we will define a class of landings
which are composed by many very good
branches, but which are controlled to an ever greater detail than excellent
landings.  Their design will allow to estimate their hyperbolicity if
truncated outside a relatively small initial segment.

We define the set of cool landings
$LC(n_0,n) \subset \Omega$, $n_0,n \in \N$,
$n \geq n_0$ as the set of all $\d=(j_1,...,j_m)$ in $LE(n_0,n)$
satisfying

\begin{description}

\item[LC1] (Starts very good) $j_i \in VG(n_0,n)$, $1 \leq i \leq
c_{n-1}^{-a^2/2}$.

\comm{
\item[LC2] (Short times are sparse in large enough initial segments)
For $c^{-a/2}_{n-1} \leq k \leq m$
$$
\#\{1 \leq i \leq k,\, r_n(j_i)<c^{-a/2}_{n-1}\}<
(6 \cdot 2^n) c^{a/3}_{n-1} k,
$$
}

\item[LC2] (Not very good moments are sparse in large enough initial
segments)
For all $c^{-a^2/4}_{n-1}<k \leq m$
$$
\#\{1 \leq i \leq k,\, r_n(j_i)<c^{-a/2}_{n-1}\} <
(6 \cdot 2^n) c^{a/3}_{n-1} k,
$$

\comm{  
\item[LC3] (Not very good moments are sparse in large enough initial
segments)
For all $c^{-a^2/4}_{n-1}<k \leq m$
$$
\#\{i|j_i \notin VG(n_0,n)\} \cap \{1,...,k\}<
(6 \cdot 2^n) c^{a^2}_{n-1} k,
$$
}

\item[LC3] (Bad moments are sparse in large enough initial segments)
For $c_{n-1}^{-n/3} \leq k \leq m$
$$
\#\{1 \leq i \leq k,\, j_i \in B(n_0,n)\} <
(6 \cdot 2^n) c^{n/6}_{n-1} k,
$$
  
\item[LC4] (Starts with no bad moments) $j_i \notin B(n_0,n)$,
$1 \leq i \leq c_{n-1}^{-n/2}$.

\end{description}

As in Lemma 7.7 of
\cite {AM}, cool landings are frequent and we get the following
parameter estimate analogous to Lemma 7.8 of \cite {AM}.  The ideas of this
estimate are quite similar to the case of standard landings.

\begin{lemma} \label {lands cool}

With total probability, for all
$n_0$ big enough, for all $n$ big enough we
have $R_n(0) \in LC(n_0,n)$.

\end{lemma}

Let us now show that cool landings inherit hyperbolicity from very good
returns.  This result corresponds to Lemma 7.11 of \cite {AM}, but the
proof of this fact needs adjustments for
big quasisymmetric constants, so we provide it here.

\begin{lemma} \label {cool hyperbolicity}

With total probability, if $n_0$ is sufficiently big, for all $n$
sufficiently big, if $\d \in LC(n_0,n)$ then for all
$c_{n-1}^{-4/(n-1)}<k \leq l_n(\d)$,
$$
\inf_{C^\d_n} \frac {\ln |Df^k|} {k} \geq
\frac {\lambda_{n_0}} {2}.
$$

\end{lemma}

\begin{pf}

Fix such $\d \in LC(n_0,n)$, and let $\d=(j_1,...,j_m)$.

Let
$$
a_k=\inf_{C^\d_n} \frac {\ln |Df^k|} {k}.
$$

Analogously to
Lemma \ref {very good partial time}, we define $m_k$ as the
number of full returns before $k$, that is,
the biggest integer such that
$$
\sum_{i=1}^{m_k} r_n(j_i) \leq k.
$$
We define
$$
\beta_k=\sum_{\ntop {1 \leq i \leq m_k,}
{j_i \in VG(n_0,n+1)}} r_n(j_i),
$$
(counting the time up to $k$ spent in complete very good returns)
and
$$
i_k=k-\sum_{i=1}^{m_k} r_n(j_i).
$$
(counting the time in the incomplete return at $k$).

Let us then consider two cases: small $m_k$ ($m_k<c_{n-1}^{-a^2/2}$)
and otherwise.

{\it Case 1 ($m_k<c_{n-1}^{-a^2/2}$).}
The idea of the first case is that all full returns are very good by LC1,
and the incomplete time is also part of a very good return.

Since full very good returns are very hyperbolic by G1 and
very good returns are good, we
just have to worry about possibly losing hyperbolicity in the incomplete
time.  To control this, we introduce the queue (or tail)
$q_k=\inf_{C^\d_n} \ln |Df^{i_k} \circ f^{k-i_k}|$.
We have $-q_k<-\ln(c_{n-1}^{1/3}
c_{n-1}^3)$ by VG and
Lemma \ref {lower bound}.  Let us split again in two cases:
$i_k$ big or otherwise.

{\it Subcase 1a ($i_k>c_{n-1}^{-4/(n-1)}$).}
If the incomplete time is big, we can use
G2 to estimate the hyperbolicity of the incomplete time (which is part of
a very good return).  The reader can easily
check the estimate in this case.

{\it Subcase 1b ($i_k<c_{n-1}^{-4/(n-1)}$).}
If the incomplete time is not big, we can not use G2 to estimate $q_k$,
but
in this case $i_k$ is much less than $k$: since $k>c_{n-1}^{-4/(n-1)}$, at
least one return was completed ($m_k \geq 1$), and since it must be very
good we conclude that $k>c_{n-1}^{-a/2}$ by LS1, so
$$
a_k>\lambda_{n_0} \frac {(1+2^{n_0-n})} {2} \cdot \frac {k-i_k} {k}-
\frac {-q_k} {k}>\frac {\lambda_{n_0}} {2}.
$$

{\it Case 2 ($m_k>c_{n-1}^{-a^2/2}$).}
For an incomplete time we still have $-q_k<-\ln(c_n c_{n-1}^3)$,
so $-q_k/k<c_{n-1}^{a^2/3}$.

\comm{
In this case,
$$
a_k>\lambda_{n_0} \frac {1+2^{n_0-n}} {2} \cdot
\frac {\beta_k} {k} - \frac {-q_k} {k},
$$
since very good returns are good and even very good returns have derivative
at least $1$, so we just have to show that $1-\beta_k/k$ is very small.
}

Arguing as in Lemma \ref {very good partial time}, we split
$k-\beta_k-i_k$
(time of full returns which are not very good) in part relative to
bad returns $h_k$ and in part relative to returns that are
not very good or bad (which must be fast) $l_k$.  Using LC3 and LC4
to bound the number of bad returns and LS2 to bound their time,
we get
$$
h_k<c_{n-1}^{-3b} c_{n-1}^{n/7} m_k,
$$
and using LC1 and LC2 we have
$$
l_k<c_{n-1}^{-a/2} c_{n-2}^{-4b} (6 \cdot 2^n) c_{n-1}^{a^2} m_k,
$$

By LC1 and LC2 again, using LS1 to estimate
the time of a very good return by $c_{n-1}^{-a/2}$, we have
that $\beta_k>c_{n-1}^{-a/2} m_k/2$, thus we get
\begin{equation} \label {hklkbetak}
\frac {h_k+l_k} {\beta_k}<c_{n-1}^{a^2/2},
\end{equation}
which is very small.

On the other hand,
$\beta_k>c_{n-1}^{-a/2}c_{n-1}^{-a^2/2}/2$ by hypothesis on $m_k$.
Let us split in three cases according to the behavior of $i_k$.

{\it Subcase 2a ($i_k$ not very good or bad).}
In this case, $i_k<c_{n-1}^{-a/2}c_{n-2}^{-4b}$, so
$i_k/\beta_k$ is very small, and we actually have
$1-\beta_k/k<c_{n-1}^{a^2/10}$.
Since very good returns are good and even
not very good returns have derivative
at least $1$,
\begin{equation} \label {2c}
a_k>\lambda_{n_0} \frac {1+2^{n_0-n}} {2} \cdot
\frac {\beta_k} {k} - \frac {-q_k} {k}>\frac {\lambda_{n_0}} {2}.
\end{equation}

{\it Subcase 2b ($i_k$ very good).}
If $i_k$ is very good and $i_k>c_{n-1}^{-4/(n-1)}$, we can reason
as in Subcase 1a that $G2$ can be used
for the estimate of $q_k$ so that we have
$$
a_k>\lambda_{n_0} \frac {1+2^{n_0-n}} {2} \cdot
\frac {\beta_k} {k}+\frac {i_k} {k} \cdot \frac {\lambda_{n_0}} {2}<\frac
{\lambda_{n_0}} {2}
$$
by (\ref {hklkbetak}).

If $i_k \leq c_{n-1}^{-4/(n-1)}$, then $i_k/\beta_k$ is very small and so
$1-\beta_k/k<c_{n-1}^{a^2/10}$, and we obtain (as in Subcase 2a)
estimate (\ref {2c}).

{\it Subcase 2c ($i_k$ bad).}
If $i_k$ is bad, by LC4 we have that $m_k>c_{n-1}^{-n/2}$,
but $i_k<c_{n-1}^{-3b}$ by LS2, so $i_k/\beta_k$ is
very small again and we have
$1-\beta_k/k<c_{n-1}^{a^2/10}$,
so estimate (\ref {2c}) applies and we are done.
\end{pf}

\subsection{Collet-Eckmann}

Since the critical point always falls in cool landings
(see Lemma \ref {lands cool}), the Collet-Eckmann condition
follows easily from Lemma \ref {cool hyperbolicity}
(which guarantees gain of derivative after large truncations), 
together with Lemma \ref {lower bound}, which controls loss of derivative
at small truncations.  This argument is identical to the one in
\S 8.1 of \cite {AM}, but we reproduce it here for the convenience of the
reader.

Let
$$
a_k=\frac {\ln |Df^k(f(0)))} {k}
$$
and $e_n=a_{v_n-1}$.

It is easy to see that if $n_0$ is big enough such that both Lemmas \ref
{lands cool} and \ref {cool hyperbolicity} we obtain for $n$ big enough
that
$$
e_{n+1} \geq e_n \frac {v_n-1} {v_{n+1}-1}+\frac {\lambda_{n_0}} {2} \cdot
\frac {v_{n+1}-v_n} {v_{n+1}-1}
$$
and so
\begin{equation} \label {hy1}
\liminf_{n \to \infty} e_n \geq \frac {\lambda_{n_0}} {2}.
\end{equation}

Let now $v_n-1<k<v_{n+1}-1$.  Define
$q_k=\ln |Df^{k-v_n}(f^{v_n}(0))|$.

Assume first that $k<v_n+c_{n-1}^{-4/(n-1)}$.
From LC1 we know that $\tau_n$
is very good, so by LS1 we have
$r_n(\tau_n)>c_{n-1}^{-a/2}$, so $k$ is in the middle of this branch (that
is, $v_n \leq k \leq v_n+r_n(\tau_n)-1$).
Using that $|R_n(0)|>|I_n|/2^n$ (by Lemma \ref {away from the boundary}), we
get by Lemma \ref {lower bound}
that $-q_k<-\ln(2^{-n} c_{n-1} c_{n-1}^5)$.
We then get from $v_n>c_{n-1}^{-a}$ that
\begin{equation} \label {hy2}
a_k \geq e_n \frac {v_n-1} {k}-\frac {-q_k} {k}>\left (1-\frac {1} {2^n}
\right )e_n-\frac {1} {2^n}.
\end{equation}

If $k>v_n+c_{n-1}^{-4/(n-1)}$ using Lemma \ref {cool hyperbolicity}
we get
\begin{equation} \label {hy3}
a_k \geq e_n \frac {v_n-1} {k}+\frac {\lambda_{n_0}} {2} \cdot
\frac {k-v_n+1} {k}.
\end{equation}

Estimates (\ref {hy1}), (\ref {hy2}), and (\ref {hy3})
imply that $\liminf a_k \geq \lambda_{n_0}/2$ and so
$f$ is Collet-Eckmann.

\subsection{Recurrence}

To show that the critical point is polynomially recurrent, we can
follow the same lines from \cite {AM}.  First we look at the
essentially Markov process $R_n|(I_n \setminus I_{n+1})$,
which shows that with total probability, most
(in the $\g$-qs sense) points in $I_n$ approach $0$ with a polynomial rate
(the exponent must be chosen according to $\g$) until the
first time they fall in $I_{n+1}$.  More precisely, we show (after
transferring to the parameter) the following estimate (analogous to
Corollary 8.3 of \cite {AM}).

\begin{lemma} \label {distance}

With total probability, for $n$ big enough and for
$1 \leq i \leq s_n$,
$$
\frac {\ln |R_n^i(0)|} {\ln(c_{n-1})} <
b^2 \left( 1+\frac {\ln(i)} {\ln(c_{n-1}^{-1})} \right).
$$

\end{lemma}

To obtain the polynomial recurrence for $f$ we relate the return times in
terms of $R_n$ to return times in terms of $f$.  In other words, letting 
$k_i$ be such that $R_n^i(0)=f^{k_i}(0)$, we must relate $k_i$ and $i$.  
It is enough to do the estimate for a cool landing and we obtain the
following estimate (as in
Corollary 8.5 of \cite {AM}).

\begin{lemma} \label {high time}

With total probability, for $n$ big enough and for
$1 \leq i \leq s_n$,
$$
\frac {\ln(k_i)} {\ln(c_{n-1}^{-1})} >
a/3 \left( 1+\frac {\ln(i)} {\ln(c_{n-1}^{-1})} \right).
$$

\end{lemma}

Let now $v_n \leq k<v_{n+1}$.  If $|f^k(0)|<k^{-3b^3}$ we have
$f^k(0) \in I_n$ and so $k=k_i$ for some $i$.
It follows from Lemmas \ref {distance} and \ref {high time} that
$$
|f^{k_i}(0)|>k^{-3b^3}_i.
$$
This concludes the proof of polynomial recurrence.  We notice that
polynomial lower bounds are easily obtained: considering
$|R_n(0)|=|f^{v_n}(0)|<c_{n-1}$ and using $v_n<c_{n-1}^{-b}$
we get
$$
\limsup_{n \to \infty} \frac {-\ln |f^n(0)|} {\ln n} \geq a.
$$

\end{document}